\documentclass[a4paper, 12pt]{article}
\usepackage{enumerate, theorem}
\usepackage{amsmath, amsfonts, amssymb}
\usepackage[height=22.5cm, width=15cm]{geometry}
\usepackage[all]{xy}
\usepackage{mathrsfs, yfonts}

\DeclareMathOperator{\id}{id}

\DeclareMathOperator{\Sym}{Sym}

 \def\C{{\mathbb C}}

\newcommand{\parag}[1]{\paragraph{\sc{#1.}} }

\newtheorem{thm}{Theorem}[subsection]
\newtheorem{defn}[thm]{Definition}
\newtheorem{cor}[thm]{Corollaire}
\newtheorem{prop}[thm]{Proposition}
\newtheorem{lemma}[thm]{Lemma}

\begin{document}

\date{07/04/15}

\author{Daniel Barlet\footnote{Institut Elie Cartan : Alg\`{e}bre et G\'eom\`{e}trie,  \newline
Universit\'e de Lorraine, CNRS UMR 7502  and  Institut Universitaire de France.}.}

\title{Strongly quasi-proper maps and the f-flattning theorem.}

\maketitle

\parag{Abstract} We complete and precise the results of [B.13] and we prove a strong version of the semi-proper direct image theorem with values in the space $\mathcal{C}_{n}^{f}(M)$ of finite type closed $n-$cycles in a complex space $M$. We describe  the strongly quasi-proper maps as the class of holomorphic surjective maps which admit a meromorphic family of fibers and  we prove stability properties of this class. In the Appendix we give a direct and short proof of D. Mathieu's flattning theorem (see [M.00]) for a strongly quasi-proper map which is easier and more accessible.\\

\parag{AMS Classification 2010} 32-H-02, 32-H-04, 32-C-25, 32-H-quasi-35, 32-15, 32-K-05, 32-L-05.

\parag{Key Words} Finite Type Cycles, Geometrically f-Flat map, Strongly Quasi-Proper map, Stein Factorization, Semi-Proper Direct Image. \newpage

\tableofcontents

\section{Introduction.}
It is now a classical result that there exists a geometric flattening theorem for a proper surjective holomorphic map between irreducible complex spaces. This is a rather simple consequence of the existence of the cycle's space of a complex space (see [B.75] and [B-M] ch.IV section 9). This fact can be seen as the meromorphie along the subset of ``big fibers'' of the ``fiber map'' for such a morphism. Of course this is a geometric version of the deep flattening theorem of H. Hironaka (see [H.75]). We consider an analogous statement for a large class of surjective morphisms between irreducible complex spaces with non compact fibers: the strongly quasi-proper maps. The main problem here is the fact that the standard notion of a quasi-proper map is not strong enough for our purpose : the strict transform of a quasi-proper map by a (proper) modification of the target space is no longer quasi-proper in general (see the example following the proposition \ref{SQP basic}). This comes from the fact that the quasi-properness is not sufficient to control the limits of generic fibers near a`` big fiber''. This motivate to work with the notion of  strongly proper maps introduced in [M.00], [B.08]  in an implicit way and in [B.13] explicitly.\\
We  complete and precise  the study of the class ``strongly quasi-proper maps'' introduced in {\it loc. cit.} which enjoys several interesting stability properties and which is an useful tool to prove the existence of meromorphic quotients in the category of reduced complex spaces (see [B.13]).\\
The main results obtained here are :
\begin{enumerate}
\item A stronger criterium for proving the SQP property and the fact that this property  is equivalent to the existence of a  meromorphic fiber map, see the proposition \ref{SPQ dense} and the theorem \ref{David's}.
\item A more precise version of a the semi-proper direct image with value in the space $\mathcal{C}_{n}^{f}(M)$, see the theorem \ref{semi-proper direct image bis} which uses a new analytic continuation result given in the theorem \ref{cycles}.
\item The strong stability theorem for SQP maps, see the theorem \ref{stab. tot}.
\end{enumerate}

\section{Analytic structure on $\mathcal{C}_{n}^{f}(M)$.}

\subsection{Definitions and examples of analytic subsets.}

We suppose that the notion of analytic family of $n-$cycles is known and also the definition of the topology of the space $\mathcal{C}_{n}^{loc}(M)$ for a complex space $M$ (see [B.M] ch. IV section 2). \\
Recall that a $n-$cycle has finite type if it has only finitely many irreducible components. The subset of finite type $n-$cycles in $M$ is denoted  $\mathcal{C}_{n}^{f}(M)$. We define the topology on it as follows :\\
For $W := (W_{1}, \dots, W_{m})$ a finite set of relatively compact open sets in $M$ let $\Omega(W)$  the subset of $n-$cycles $C$  in $M$ such that each irreducible component of $\vert C\vert$ meets each \ $W_{i}$ for $i = 1, \dots, m$. Remark that finite intersection of some $\Omega(W)$ is again of the form $\Omega(W)$. Now define open sets in $\mathcal{C}^{f}_{n}(M)$ as union of subsets of the form  \ $\mathcal{U}\cap \Omega(W)$ \ where $\mathcal{U}$ is an open set in $\mathcal{C}^{loc}_{n}(M)$.\\
So the inclusion $\mathcal{C}_{n}^{f}(M) \to \mathcal{C}_{n}^{loc}(M)$ is continuous, but it is not a homeomorphism on its image (with the induced topology). For instance, it is easy to see that a continuous family of $n-$cycles $(X_{s})_{s \in S}$ in $M$ parametrized by a Hausdorff topological space $S$, such that each cycle is of finite type is f-continuous (so corresponds to a continuous map $S \to \mathcal{C}_{n}^{f}(M)$ with the topology defined above) if and only if its set-theoretic graph $\vert G\vert \subset S \times M$ is quasi-proper over $S$ ; here we use the following definition of quasi-proper, valid as long as we know that the fibers of the continuous map  $\pi : G \to S$ are complex analytic subsets (here $G \subset S\times M$ is closed  and the fibers are cycles in the complex space $M$) :\\
For any point $s_{0}\in S$ there exists a neighborhood $S_{0}$ of $s_{0}$ in $S$ and a closed $S_{0}-$proper  subset $K$ in $\pi^{-1}(S_{0})$ such that for any $s \in S_{0}$ any irreducible component of the fiber $\pi^{-1}(s)$ meets $K$.

\begin{lemma}\label{Denom.}
For any complex space $M$ and any integer $n$ the topology of the space $\mathcal{C}_{n}^{f}(M)$ has a countable basis\footnote{By definition a complex space is countable at infinity.}.
\end{lemma}

So this space is metrizable and then any point has a countable basis of closed neighborhoods. 

\parag{Proof} This an easy consequence of the analogous result for the topology of $\mathcal{C}_{n}^{loc}(M)$ which is proved in [B-M] IV section 2.4: \\
Let $D$ be a countable dense subset in $M$ and consider also a countable covering of $M$ by domains of charts $(U_{p}, j_{p})_{p \in \mathbb{N}}$. Now for each point $d \in D$ let $p(d)$ be the smallest integer $p$  such that $d \in U_{p}$. Define for $r \in \mathbb{Q}^{+*}$ small enough the pull-back $P_{d,r}$ in $U_{p(d)}$ of the polydisc with center $j_{p(d)}(d)$ and radius $r$. Then take a countable basis $\mathcal{U}_{\nu}, \nu \in \mathbb{N}$,  for the topology of $\mathcal{C}^{loc}_{n}(M)$ and consider the  family $(\Omega(W_{j})\cap \mathcal{U}_{\nu})_{j,\nu}$ where the countable family of relatively compact open sets $W_{j} =: (W_{j_{1}}\dots W_{j_{\mu}})$ runs over the finite subsets of the countable family of the $P_{d,r} \subset\subset M$,  gives a countable basis for the topology of $\mathcal{C}_{n}^{f}(M)$.$\hfill \blacksquare$\\

\begin{prop}\label{f-conv.}
Let \ $M$ \ be a reduced complex space and consider a sequence of finite type \ $n-$cycles \ $(C_{m})_{m \geq 0}$ \ in \ $M$ \ with the following properties :
\begin{enumerate}[i)]
\item there exists a compact set \ $K$ \ in \ $M$ \ such that each irreducible component of each \ $C_{m}$ \ meets \ $K$.
\item The sequence \ $(C_{m})_{m \geq 0}$ \ converges for the topology of \ $\mathcal{C}_{n}^{loc}(M)$ \ to a \ $n-$cycle \ $C$.
\item The cycle \ $C$ \ is in \ $\mathcal{C}_{n}^{f}(M)$.
\end{enumerate}
Then the sequence \ $(C_{m})_{m \geq 0}$ \ converges to \ $C$ \ in the topology of \ $\mathcal{C}_{n}^{f}(M)$.
\end{prop}


\parag{Proof} We begin by proving the case where each $C_{m}$ is irreducible. Note first that \ $C$ \ is not the empty \ $n-$cycle because, up to pass to a subsequence, we may assume that for each \ $m \geq 0$ \ we can choose a point in \ $x_{m} \in C_{m} \cap K$ \ in order that the sequence \ $(x_{m})_{m \geq 0}$ \ converges to a point \ $x \in K$. Then we have \ $x \in \vert C\vert$.\\
Let \ $W := (W_{1}, \dots, W_{\mu})$ be a finite collection of relatively compact open sets open \ $M$ \ such that each irreducible component of \ $\vert C\vert$ meets each $W_{i}$ for $i = 1, \dots, \mu$. Choose also $F_{j}, j \in J$,  a finite number of $n-$scales on $M$ adapted to $C$ and put $l_{j} := \deg_{F_{j}}(C)$. Define the open set of $\mathcal{C}^{f}_{n}(M)$ containing $C$
$$ \mathcal{V} := \Omega(W_{1}, \dots, W_{\mu}) \cap_{j \in J} \Omega_{l_{j}}(F_{j}) .$$
 Choose also for each $i \in [1,\mu]$ a  $n-$scales \ $E_{i} :=(U_{i}, B_{i}, j_{i})$ \   on \ $W_{i}$ \ adapted to \ $C$ \ and such that \ $\deg_{E_{i}}(C) =k_{i} \geq 1 $. This is possible because \ $C$ \ has at least one irreducible component  and such an irreducible component meets each \ $W_{i}$. So for \ $m \geq m_{0}$ \ the scales $F_{j}$ and  \ $E_{i}$ \ will be adapted to \ $C_{m} $ \ and we shall have  for all $j \in J$ and all $i \in [1, \mu]$ 
$$l_{j} = \deg_{F_{j}}(C) = \deg_{F_{j}}(C_{m}) \quad {\rm and} \quad        \deg_{E_{i}}(C_{m}) = \deg_{E_{i}}(C) = k_{i} \geq 1.$$
 This shows that for \ $m \geq m_{0}$ \  the unique  irreducible component of \ $C_{m}$ meets each  $W_{i}$. So, for $m \geq m_{0}$ each $C_{m}$ for $m \geq m_{0}$ is in the given open  set $\mathcal{V}$ of $\mathcal{C}^{f}_{n}(M)$. So  the sequence \ $(C_{m})_{m \geq 0}$ \ converges to \ $C$ \ in the sense of  the topology of \ $\mathcal{C}_{n}^{f}(M)$. This is enough to conclude the proof in this case.\\
 Consider now the general case. We shall prove first that it is enough to prove the following claim :
\parag{Claim} For any relatively compact open set $W$ meeting each irreducible component of $C$, there exists $m_{1}$ such for $m \geq m_{1}$ any irreducible component of any $C_{m}$ meets $W$.
\parag{The claim is enough} Consider as before an open set
 $$\mathcal{V} := \Omega(W_{1}, \dots, W_{\mu}) \cap_{j \in J} \Omega_{l_{j}}(F_{j})$$
   in $\mathcal{C}^{f}_{n}(M)$  containing $C$. If the claim is true we find an integer $m_{2}$ such that for each $m \geq m_{2}$ we have $C_{m} \in \Omega(W_{1}, \dots, W_{\mu})$. Now the convergence of the sequence for the topology of $\mathcal{C}_{n}^{loc}(M)$ gives an integer $m_{3}$ such that for each $m \geq m_{3}$ we have $C_{m} \in \cap_{j \in J} \Omega_{l_{j}}(F_{j})$. So the convergence for the topology of $\mathcal{C}^{f}_{n}(M)$ is proved.
\parag{proof of the claim} Take any open set \ $W \subset M$ \ such that any irreducible component of \ $C$ \ meets \ $W$.  If there is infinitely many $m \geq 0$ for which an irreducible component $\Gamma_{m}$ of $C_{m}$ does not meet $W$ then extract a sub-sequence $(\Gamma_{p})$ of irreducible components in  this sequence which converges  in the topology of \ $\mathcal{C}_{n}^{loc}(M)$ \ to a {\bf non empty}\footnote{remember that each $\Gamma_{m}$ meets the compact set $K$.} cycle \ $\Gamma \leq C$. But any irreducible component of \ $\Gamma$ \ has to meet \ $W$ \ giving a contradiction as \ $\Gamma_{p}\cap W = \emptyset$ \ by construction. So the claim is proved and the proof is complete.$\hfill \blacksquare$\\

\begin{cor}\label{compact}
Let \ $M$ \ be a reduced complexe space and let \ $A \subset \mathcal{C}_{n}^{f}(M) \setminus \{\emptyset\}$ \ be a  compact subset in \ $\mathcal{C}_{n}^{loc}(M)$. Assume that the following condition is fullfilled :
\begin{itemize}
\item There exists a compact set \ $K \subset M$ \ such for any irreducible component \ $\Gamma$ \ of each  \ $C \in A$ \ the intersection \ $\Gamma \cap K$ \ is not empty.$\hfill (@)$
\end{itemize}
Then  \ $A$ \ is  compact in \ $\mathcal{C}_{n}^{f}(M)$.\\
\end{cor}

Note that conversely any compact subset in $\mathcal{C}_{n}^{f}(M) \setminus \{\emptyset\}$ is compact in $\mathcal{C}_{n}^{loc}(M)$ and satisfies $(@)$ : choose for each $C \in A$ an open set $\mathcal{U}_{C}\cap \Omega(W_{C})$ containig $C$. As $A$ is compact we can extract a finite sub-covering given by $C_{1}, ..., C_{N}$ and take $K := \overline{\cup_{j=1}^{N} \ W_{C_{j}}}$.

\parag{Proof} Note that any limit of a sequence in \ $A$ \ for the topology of \ $\mathcal{C}_{n}^{loc}(M)$ \ is not the empty cycle, as it contains a point in \ $K$.  Now, by assumption, for any sequence in \ $A$ \ we may find a sub-sequence converging for the topology of \ $\mathcal{C}_{n}^{loc}(M)$ \  to a cycle in \ $\mathcal{C}_{n}^{f}(M)$. So, applying the proposition \ref{f-conv.}, we conclude that such a sub-sequence converges for the topology of \ $\mathcal{C}_{n}^{f}(M)$.$\hfill \blacksquare$

\begin{defn}\label{hol in}
Let $S$ be a reduced complex space. A map $\varphi : S \to \mathcal{C}_{n}^{f}(M)$ will be called {\bf holomorphic} when it classifies a f-analytic family of $n-$cycles in $M$.
\end{defn}

Recall that an analytic family of $n-$cycles in $M$ is f-analytic if and only if its set theoretic graph $\vert G \vert \subset S \times M$ is {\bf quasi-proper} on $S$. Of course this condition implies that each cycle in the family is a finite type cycle. But the quasi-properness of the graph asks more : on any compact set $L$ in $S$ we can find a compact set $K$ in $M$ such any irreducible component of any cycle parametrized by a point $s \in L$ has to meet $K$.

\begin{defn}\label{hol to} Let \ $\mathcal{U}$ be an open set in $\mathcal{C}_{n}^{f}(M)$.
A continuous map  $ f : \mathcal{U} \to T$  to a Banach analytic set \  $ T$  is called {\bf holomorphic} if and only if for any reduced complex space $S$ and any  holomorphic map $\varphi : S \to \mathcal{U}$ the composed map
 $$ f\circ\varphi : S \to T$$
 is holomorphic. We shall say that a closed subset $\mathcal{X}$  in an open set \ $\mathcal{U}$ of $\mathcal{C}_{n}^{f}(M)$ is {\bf analytic} when there exists locally on \ $\mathcal{U}$ an holomorphic map $f$ with values in a Banach space $E$  such that $\mathcal{X} = f^{-1}(0)$. \\
A continuous map of an analytic subset $\mathcal{X}$ of \  $\mathcal{U}$  to a Banach space  $ E$ will be called {\bf holomorphic} if it is locally induced on $\mathcal{X}$ by a holomorphic map of $\mathcal{C}_{n}^{f}(M)$  with values in  $ E$.
\end{defn}

The next lemma gives a first example of a closed analytic subset in $\mathcal{C}_{n}^{f}(M)$.

\begin{lemma}\label{reduced}
Let $NR := \{C \in \mathcal{C}_{n}^{f}(M)\ / \  C \not= \vert C \vert \}$ be the subset of non reduced cycles. Then $NR$ is a closed analytic subset in $\mathcal{C}_{n}^{f}(M)$.
\end{lemma}

\parag{Proof} As the empty $n-$cycle is open and closed in $\mathcal{C}_{n}^{f}(M)$ it is enough to consider non empty cycles. Let $C_{0}$ be any non empty  cycle in $\mathcal{C}_{n}^{f}(M)$. Choose for each irreducible component $\Gamma_{i}$ of $C_{0}$ an $n-$scale $E_{i} := (U_{i}, B_{i},j_{i})$ on $M$ adapted to $C_{0}$ such that the degree of $\vert C_{0}\vert$ and  $\Gamma_{i}$ in $E_{i}$ are equal to $1$, and note $k_{i} := \deg_{E_{i}}(C_{0})$. Remark that $C_{0}$ is reduced if and only if we have $k_{i} = 1$ for each $i \in I$. Let $W := \cup_{I \in I} \ j_{i}^{-1}(U_{i}\times B_{i})$ and $\mathcal{V} := \Omega(W) \cap \big(\cap_{i\in I} \ \Omega_{k_{i}}(E_{i})\big)$. Then a cycle $C \in \mathcal{V}$ is not reduced if and only if there exists at least one  $i \in I$ such that $C \cap  j_{i}^{-1}(U_{i}\times B_{i})$ is not reduced. As each map $\mathcal{V} \to H(\bar U_{i}, \Sym^{k_{i}}(B_{i}))$ is holomorphic, the proof is consequence of the following claim :
\parag{Claim} \begin{itemize}
\item The subset of $H(\bar U, \Sym^{k}(B))$ corresponding to non reduced cycles in $U \times B$ is a closed analytic subset.
\end{itemize}
Consider the discriminant map $\Delta_{0} : \Sym^{k}(\C^{p}) \to S^{k.(k-1)}(\C^{p})$ defined by
 $$(x_{1}, \dots, x_{k}) \mapsto \prod_{1\leq i < j \leq k} \quad (x_{i} - x_{j})^{2} .$$
 It is induced by a polynomial map $\oplus_{i=1}^{k} \ S^{i}(\C^{p}) \to S^{k.(k-1)}(\C^{p})$ thanks to the symmetric function theorem (for this vector case see [B-M] theorem II 4.2.7), and so we have a holomorphic map $ \Delta : H(\bar U, \Sym^{k}(B)) \to H(\bar U, S^{k.(k-1)}(\C^{p}))$ given by $f \mapsto \Delta_{0}\circ f$. Of course, if $f \in H(\bar U, \Sym^{k}(B))$ defines a non reduced cycle in $U\times B$ we have $\Delta_{0}\circ f = 0$ in $H(\bar U, S^{k.(k-1)}(\C^{p}))$. Conversely, if $\Delta_{0}\circ f = 0$ in $H(\bar U, S^{k.(k-1)}(\C^{p}))$, then choose a point $t_{0}$ which is not a ramification point for the reduced cycle which is the support of the cycle $X$  associated to $f$. Then locally $\vert X \vert$ is the disjoint union of $h \leq k$ graphs of holomorphic functions $f_{1}, \dots, f_{h} : V(t_{0}) \to B$ which are two by two disjoint. If $X$ is reduced, we have $h = k$ and on the open neighborhood $V(t_{0})$ we have
 $$ \Delta_{0}\circ f = \prod_{1\leq i < j \leq k} \quad (f_{i} - f_{j})^{2} $$
 which nowhere vanishes on $V(t_{0})$. Contradiction. So the claim is proved.$\hfill \blacksquare$\\

The second interesting example of an analytic subset of $\mathcal{C}_{n}^{f}(M)$ and of holomorphic map on it  is given by the following proposition.

\begin{prop}\label{relatif}
Let $\pi : M \to S$ be a holomorphic map between two reduced complex spaces. Define $\mathcal{C}_{n}^{f}(\pi)$ as the subset of $\mathcal{C}_{n}^{f}(M)$ of the $n-$cycles contained in a fiber of $\pi$. Then $\mathcal{C}_{n}^{f}(\pi)$ is an analytic (closed) subset in $\mathcal{C}_{n}^{f}(M)$ and the obvious map $p : \mathcal{C}_{n}^{f}(\pi) \to S$ is holomorphic.
\end{prop}

\parag{Proof} First we shall show that the complement of $\mathcal{C}_{n}^{f}(\pi)$ is open: \\
take $C_{0}\not\in \mathcal{C}_{n}^{f}(\pi)$. Then it contains two points $x$ and $y$ such that $\pi(x) \not= \pi(y)$. Take two scales $E := (U,B,j)$ and $E' := (U',B',j')$  adapted to $C_{0}$ such that $x$ is in the center $j^{-1}(U\times B)$ of $E$ and $y$ in the center $(j')^{-1}(U'\times B')$ of $E'$, small enough such that $\pi(j^{-1}(\bar U\times \bar B))$ and $\pi((j')^{-1}(\bar U'\times \bar B') $ are disjoints.  Note that the degrees of $C_{0}$ in $E$ and in $E'$ are positive. Now a cycle $C$ such that $E$ and $E'$ are adapted to $C$ with the same degrees as $C_{0}$ in these scales cannot be in $\mathcal{C}_{n}^{f}(\pi)$ because it has two points with different images by $\pi$. This defines an open set in the complement of $\mathcal{C}_{n}^{f}(\pi)$ containing $C_{0}$.

To obtain a local holomorphic equation for $\mathcal{C}_{n}^{f}(\pi)$, recall the following facts :
\begin{enumerate}[i)]
\item For any $n-$scale  $E := (U,B,j)$ on $M$  let $\Omega_{k}(E)$ be the open set in $\mathcal{C}_{n}^{loc}(M)$ of cycles for which $E$ is adapted and the degree in $E$ is $k$. Then the map
$$ r_{E,k} : \Omega_{k}(E) \to H(\bar U, \Sym^{k}(B)) $$
is holomorphic. This is a consequence of the definition of an analytic family of cycles !
\item If for a given cycle $C_{0}$ we have adapted scales $E_{1}, \dots, E_{m}$ such that any irreducible component of $C_{0}$ meets the union $W$  of the centers of the $E_{i}, i \in [1,m]$; then the subset of  \ $\big[\cap_{i \in [1,m]} \Omega_{k_{i}}(E_{i})\big] \cap \Omega(W)$, where $k_{i} := deg_{E_{i}}(C_{0})$, is an open set in $\mathcal{C}_{n}^{f}(M)$  and the holomorphic map  $\prod_{i \in [1,m]} r_{E_{i},k_{i}}$ is injective on this open set.
\end{enumerate}
 Then the following lemma allows to conclude. $\hfill \blacksquare$

\begin{lemma}\label{rel.}
Let $U$ and $B$ be relatively compact polydiscs respectively in $\C^{n}$ and $\C^{p}$, and let $\pi : W \to F$ be a holomorphic map of an open neighborhood $W$ of \ $\bar U \times \bar B$ to a Banach space $F$. Then the subset $\mathcal{X}$ of $H(\bar U, \Sym^{k}(B))$ of multiform graphs contained in a fiber of $\pi$ is a closed Banach analytic subset of $H(\bar U, \Sym^{k}(B))$. Moreover, the map $\varphi : \mathcal{X} \to F$ defined by sending  $X \in \mathcal{X}$  to the point $\pi(X) \in F$ such that $\vert X\vert \subset \pi^{-1}(\pi(X))$ is holomorphic.
\end{lemma}

\parag{proof} Consider for each $h \in [1,k]$ the holomorphic map
 $$N_{h}(\pi) : \Sym^{k}(W)  \to S^{h}(F)$$
  given by the $h-$th Newton symmetric function $(z_{1}, \dots, z_{k}) \mapsto \sum_{j=1}^{k} \ \pi(z_{j})^{h}$, where $S^{h}(F)$ is the $h-$th symmetric power of $F$\footnote{That is to say the Banach space generated by $x^{h}, x \in F$, in the Banach space of continuous homogeneous polynomials  of degree $h$ on the dual $F'$ of $F$.}. These maps are holomorphic and for $f \in H(\bar U, \Sym^{k}(B))$ we can compose the associated map  $\tilde{f} \in H(\bar U, \Sym^{k}(W))$, sending $t$ to $((t,x_{1}), \dots, (t,x_{k}))$ if $f(t) := (x_{1}, \dots, x_{k})$, with $\oplus_{h=1}^{k} N_{h}(\pi)$ to obtain a holomorphic map
  $$ \Phi :  H(\bar U, \Sym^{k}(B)) \longrightarrow H(\bar U, \oplus_{h=1}^{k} S^{h}(F)) .$$
  Now fix a point $t_{0}\in U$ and a non empty open  polydisc $U' \subset\subset U$ and consider the holomorphic maps
  $$ \Psi : H(\bar U, \oplus_{h=1}^{k} S^{h}(F)) \to H(\bar U', L(\C^{n}, \oplus_{h=1}^{k} S^{h}(F))) $$
  and
  $$ \chi :  H(\bar U, \oplus_{h=1}^{k} S^{h}(F)) \to \oplus_{h=2}^{k} S^{h}(F) $$
  defined as follows :   $\Psi(f)$ is given by the derivative of $f$ on $\bar U'$ and $\chi(f)$ is given by the collection of the \  $k^{h-1}.N_{h}(f(t_{0})) - \big( N_{1}(f(t_{0}))\big)^{h} \in S^{h}(F)$ \ for $h \in [2,k]$. These maps are holomorphic and the Banach analytic subset $ Z := \Psi^{-1}(0) \cap \chi^{-1}(0)$ is the subset in $ H(\bar U, \oplus_{h=1}^{k} S^{h}(F))$ corresponding to constant maps $ \bar U \to \oplus_{j=1}^{k} S^{h}(F))$ such that the value is of the form $ k.a \oplus k.a^{2} \oplus \dots \oplus k.a^{k}$ for some $a \in F$. So $\Phi^{-1}(Z)$ is exactly the subset $\mathcal{X}$ of $H(\bar U, \Sym^{k}(B))$.\\
  To conclude the proof, it is enough to remark that the holomorphic map   $\frac{1}{k}.\big(ev_{1}\circ \Phi\big)$ induces on $\mathcal{X}$ the desired map, where $ev_{1} :  H(\bar U, \oplus_{h=1}^{k} S^{h}(F)) \to F$ is given by evaluation at $t_{0}$ of the component  on $S^{1}(F) = F$  .$\hfill \blacksquare$\\
  
  Note  that an obvious consequence of the previous proposition is the fact that if $N$ is a closed analytic subset in $M$, then $\mathcal{C}_{n}^{f}(N)$ is a closed analytic subset of $\mathcal{C}_{n}^{f}(M)$.\\
  
 An useful variant of this result is given by the following lemma.
 
 \begin{lemma}\label{1}
Let $M$ and $S$ be reduced complex spaces. Let $(X_{s})_{s \in S}$ be an f-analytic family of $n-$cycles of $M$ parametrized by $S$. Let $T$ be an analytic set in $M$. Then the subset 
$$ Z := \{s \in S\ / \  \vert X_{s}\vert  {\rm \ has \ an\ irreducible \ component \ contained\ in \ } T\}$$
is a closed analytic subset in $S$.
\end{lemma}

\parag{proof} It is enough to consider the case of an analytic family of multiform graphs in a product $U \times B$ of polydiscs, classified by an holomorphic map
 $$f : S \times U \to \Sym^{k}(B),$$
  where  $T \subset U\times B$ is the subset  given by $T := g^{-1}(0) $ for a holomorphic function   $g : U \times B \to \C$. If $\tilde{f} : S \times U \to \Sym^{k}(U\times B)$ is the holomorphic map associated to $f$ keeping the $U-$component, then consider, for $U' \subset \subset U$, the holomorphic map
   $$\Phi : S \to H(\bar U', \C) \quad {\rm given \ by} \quad s \mapsto \big(t \mapsto Nr(g)[\tilde{f}(s,t)]\big)$$
   where $Nr(g) : \Sym^{k}(U\times B) \to \C$ is the norm\footnote{So $Nr(g)[z_{1}, \dots, z_{k}] = \prod_{j=1}^{k} \ g(z_{j})$ \ for $[z_{1}, \dots, z_{k}] \in \Sym^{k}(U\times B)$.} of $g$.\\
    Then for $U'$ non empty,  the set $\Phi^{-1}(0)$ in $S$ is the corresponding $Z$. As $\Phi$ is holomorphic, this conclude the proof.$\hfill \blacksquare$\\
    
    This lemma allows to give the following corollary to the proposition \ref{relatif}

\begin{cor}\label{une comp.}
Let $T $ be a closed analytic subset of the complex space $M$. Then the subset $\mathcal{T} \subset \mathcal{C}_{n}^{f}(M)$ of cycles having an irreducible component contained in $T$ is a closed analytic subset in $\mathcal{C}_{n}^{f}(M)$.
\end{cor}

\parag{Proof} First we want to prove that the complement of $\mathcal{T}$  is open. Let $C_{0} \not\in \mathcal{T}$ and choose a point $x_{j}, j \in [1,m]$, which is not in $T$,, and choose  in each irreducible component of the cycle $C_{0}$ a  $n-$scales $E_{j}, j \in [1,m]$, on $M \setminus T$, such that $x_{j}$ is in the center of $E_{j}$ for each $j$. Then $k_{j} := \deg_{E_{j}}(C_{0})$ is positive for each $j$. Let  $W$ be  the union of the centers of the scales $E_{j}, j \in [1,m] $ and  defined the open set $\mathcal{U}$ in $\mathcal{C}_{n}^{f}(M)$  which contains $C_{0}$  by 
$$ \mathcal{U} := \Omega(W) \cap \Big[ \cap_{j \in [1,m]} \Omega_{k_{j}}(E_{j})\Big] $$
where \ $\Omega_{k}(E)$\  denotes  the open  set of cycles \ $C$ \ such \ $E$ \ is adapted to \ $C$ with \\ $\deg_{E}(C) = k$ and $ \Omega(W)$ is the set of cycles such that any irreducible component meets $W$. Then $\mathcal{U}$  does not meet $\mathcal{T}$ and so $\mathcal{T}$ is closed.\\
To prove the analyticity of $\mathcal{T}$ consider a cycle  $C_{0}$ in $\mathcal{T}$ and again choose a point $x_{j} $ in each irreducible component of $C_{0}$ and for each $j \in [1,m]$  a  $n-$scale $E_{j}$ adapted to $C_{0}$ with center containing $x_{j}$. Construct the open set \ $\mathcal{U}$ \  as above and write in this open set the Banach analytic equations given by the previous lemma in each scale $E_{j}$. The (finite) union of the analytic sets  so defined for each $j \in [1,m]$ is equal to $\mathcal{T} \cap \mathcal{U}$, concluding the proof. $\hfill \blacksquare$\\

Note that this result may not be true in the case of an analytic family of cycles  which is not f-analytic as it is shown by the following example.

\parag{Example} Let $M := D = \{z \in \C \ / \ \vert z \vert < 1 \}$, $T := \{ 1 - \frac{1}{n}, n \in \mathbb{N}, n \geq 3\}$ and consider the family of $0-$cycles in $D$ parametrized by $D$ :
$$ X_{s} := \{ 1 - \frac{s+1}{s + m}, m \in \mathbb{N}, m \geq 3 \}\cap D  \quad {\rm for} \quad s \in D .$$
We have $X_{0} = T$ and a necessary and sufficient  condition on $s \in D$ in order that $X_{s}$ meets $T$ is that there exists $m,n \in \mathbb{N} \setminus \{0, 1,2\}$ with $ \frac{1}{n} = \frac{s+1}{s+m}$. This gives that $X_{s}$ meets $T$ if and only iff $s =\frac{p}{q}$ with $p \in \mathbb{Z}, \ q \in \mathbb{N}\setminus \{0, 1\}$ and $\vert \frac{p}{q} \vert < 1$. This is a dense set in $]-1,+1[$ !  $\hfill \square$\\

We shall conclude this paragraph by an elementary but useful result on the number of irreducible components of cycles in a f-continuous family of cycles.

\begin{defn}\label{poids}
Let $X := \sum_{i \in I} \ n_{i}.X_{i}$ be a finite type $n-$cycle in a complex space $M$, where $X_{i}$ are the irreducible components of $\vert X\vert$, $n_{i}$ are positive integers and $I$ is a finite set. We define the {\bf weight of $X$} as the integer  $w(X) := \sum_{i \in I} \ n_{i}$.
\end{defn}

Of course when $X$ is a reduced cycle, the weight of $X$ is simply the number of irreducible components of $X$.

\begin{lemma}\label{loc. bd.}
Let $M$ be a complex space and $n$ an integer. The weight function $ w : \mathcal{C}_{n}^{f}(M) \to \mathbb{N}$ is lower semi-continuous on $\mathcal{C}_{n}^{f}(M)$.
\end{lemma}

\parag{Proof} Let $C_{0}$ be a non empty cycle in $\mathcal{C}_{n}^{f}(M)$ and choose on $M$, for each irreducible component $\Gamma_{i}$  of $\vert C_{0}\vert$, an adapted scale $E_{i} := (U_{i}, B_{i}, j_{i})$ such that we have $\deg_{E_{i}}(\Gamma_{i}) = \deg_{E_{i}}(\vert C_{0}\vert) = 1$ and $\deg_{E_{I}}(C_{0}) = k_{i}$. Then we have $w(C_{0}) = \sum_{i \in I} \ k_{i}$.  Let $W := \cup_{i\in I} \  j_{i}^{-1}(U_{i}\times B_{i})$ and define the open neighbourhood $\mathcal{V}$ of $C_{0}$ in $\mathcal{C}_{n}^{f}(M) $ as :
$$ \mathcal{V} := \Omega(W) \cap (\cap_{i \in I} \  \Omega_{k_{i}}(E_{i})).$$
Then for any $C \in \mathcal{V}$ we have the inequality $ 1 \leq w(C) \leq \sum_{i \in I} \ k_{i} = w(C_{0})$. This is consequence of the fact that each irreducible component of $C$ has to meet $W$, and that the degree of $C$ in the scale $E_{i}$ is equal to $k_{i}$.$\hfill \blacksquare$\\

Then, using the existence of a meromorphic geometrically flat  Stein reduction for a meromorphic geometrically f-flat map proved in [B.13] theorem 4.2.4, we  obtain the following result.

\begin{prop}\label{red. const.}
Let $M$ be a complex space, $n$ an integer and  $\varphi : N \to \mathcal{C}_{n}^{f}(M)$ be a holomorphic map where $N$ is an irreducible complex space such that for $y$ in a dense Zariski open set  the cycle $\varphi(y)$ is reduced. Then the number of irreducible components of $\varphi(y)$  is constant on a dense open set in $N$.$\hfill \blacksquare$
\end{prop}

Note that, thanks to the lemma \ref{reduced}, if there exists one $y \in N$ such that $\varphi(y)$ is a reduced cycle, then there exists a dense Zariski open set $N'$  in $N$ such for all $y \in N'$ the cycle $\varphi(y)$ is reduced.

\parag{Proof} Let $G \subset N \times M$ be the graph of the f-analytic family of $n-$cycles of $M$ classified by $\varphi$. As the projection $\pi : G \to N$ is geometrically f-flat, we can apply the theorem 4.2.4  of [B.13], see also  the theorem \ref{David's} in  section 4.  Let $\nu : \tilde{G} \to G$ the normalization of $G$, $q : \tilde{G} \to Q$ and $p : Q \to N$ the maps given by the theorem. Then $q$ is geometrically f-flat with irreducible general fibers and $p$ is proper finite and surjective. If $k $ is the generic degree of $p$, then it is easy to see, using the lemma \ref{1} with the subset $\nu^{-1}(T) \subset \tilde{G}$, where $T \subset G$ is the subset of non normal points in $G$, that there is a  dense open subset  in $N$ on which each $\varphi(y)$ has exactly $k$ irreducible components.$\hfill \blacksquare$.\\

The next lemma shows that the graph of a f-analytic family of cycles in $M$ parametrized by an irreducible complex space is again a finite type cycle, assuming that the generic cycle is reduced. It is an easy exercice left to the reader to show  that the existence of a reduced cycle in such a family is not necessary for this result.

 \begin{lemma}\label{finitude comp.}
 Let $M$ be a complex space and $N$ an irreducible complex space. Let $\varphi : N \to \mathcal{C}_{n}^{f}(M)$ be a holomorphic map such that for some $y$ in $N$ the cycle $\varphi(y)$ is reduced. Let $G \subset N \times M$ be the graph of the f-analytic family of $n-$cycles in $M$ classified by the map $\varphi$. Then $G$ has finitely many irreducible components.
 \end{lemma}
 
 \parag{Proof} Let $N'$ be the set of points in $N$ such that $\varphi(y)$ is reduced. Then the irreducibility of $N$ and  the lemma \ref{reduced} imply that $N'$ is an dense Zariski open set in $N$. Denote $p : G \to N$ the projection. Then $G' :=  p^{-1}(N')$ is a dense Zariski open set in $G$. Then the proposition \ref{red. const.} gives that the number of irreducible components of $\varphi(p(z))$ for any $z$ in a dense open set $G'' := p^{-1}(N'')$  is constant equal to a positive integer $k$, where $N''\subset N'$ is a dense open set.\\
Note $I$ the set of irreducible components of $G$. For $i \in I$ we shall denote $G_{i}$ the corresponding irreducible component.  Let $Z $ be the analytic subset in $G$ corresponding to points which are in two distinct irreducible components of $G$. It is a closed analytic subset with no interior point in $G$ and also in any irreducible component $G_{i}$ of $G$. Now for each $i$ the image of the restriction $p_{i}$ of the projection $p$ to $G_{i}$  contains an open dense subset $ N_{i}$\footnote{By the remark 2 following the definition 3.1.4 each $p_{i}$ is in fact surjective.}, as this restriction is equidimensional between irreducible complex spaces. Now the set $I$ is countable, so the intersection $N_{0} := N'' \cap (\cap_{i \in I} \ N_{i})$ is dense by Baire's theorem. Let $N_{1}$ be the dense Zariski open set in $N$ such that $\varphi(y)$ has no irreducible component contained in $Z$.Consider now a point $y$ in $N_{0}\cap N_{1}$. Then $\varphi(p(y))$ is reduced, have $k$ irreducible components, and each irreducible component $G_{i}$ contains at least one of these components; but any irreducible component of $\varphi(p(y))$ cannot be in two different $G_{i}$ because $y$ is in $N_{1}$. So the set $I$ has at most $k$ elements.$\hfill \blacksquare$\\

\subsection{An analytic extension theorem.}

The main goal of this paragraph is to prove the following analytic continuation result.

\begin{thm}\label{cycles}
Let $M$ be a complex space and  $n$ an integer. Consider a f-continuous family  $(X_{s})_{s\in S}$ of finite type  $n-$cycles in  $M$ parametrized by a reduced complex space  $S$. Fix a point $s_{0}$ in $S$ and assume that there exists an open set $M'$ in $M$ meeting any irreducible component of $\vert X_{s_{0}}\vert$  and such that the family $(X_{s}\cap M')_{s \in S}$ is analytic at  $s_{0}$. Then there exists an open neighborhood  $S_{0}$ of $s_{0}$ in $S$ such that the family  $(X_{s})_{s\in S_{0}}$ is  f-analytic.
\end{thm}

Let us make explicit the situation of the previous theorem in term of classifying maps : we have a continuous map $ \varphi : S \to \mathcal{C}_{n}^{f}(M) $ such that the composed map $r\circ \varphi$ is holomorphic at $s_{0}$, where $r : \mathcal{C}_{n}^{f}(M)  \to \mathcal{C}_{n}^{loc}(M') $ is the restriction map. Then the statement is that $\varphi$ is holomorphic on an open neighborhood  $S_{0}$ of $s_{0}$ in $S$, assuming that $M'$ meets each irreducible component of  $\vert X_{s_{0}}\vert$.\\
Remark that the map $r$ is holomorphic\footnote{meaning that for any holomorphic map $\psi : T \to \mathcal{C}_{n}^{f}(M) $ of a reduced complex space $T$ the composed map $r\circ\psi$ is holomorphic.} so that the holomorphy at $s_{0}$ of $r\circ \varphi$ is a necessary condition for the holomorphy of $\varphi$ on an open neighborhood of $s_{0}$. The theorem says that this condition is sufficient.\\

One key point in the proof of the previous theorem is given by the following analytic continuation result.

\begin{prop}\label{fonction}
Let $S$ be a reduced complex space and let  $U_{1} \subset U_{2}$ be two open  polydiscs in  $\C^{n}$. Let  $f : S \times U_{2} \to \C$ be a continuous function, holomorphic on  $\{s\}\times U_{2}$ for each fixed $s \in S$ and assume also that the restriction of $f$ to   $S \times U_{1}$ is holomorphic. Then $f$ is holomorphic on $S \times U_{2}$.
\end{prop}

\parag{Proof of the proposition} Consider first the case where $S$ is smooth. As the problem is local on $S$ it is enough to treat the case where $S$ is an open set in some $\C^{m}$. Fix then a relatively compact open polydisc $P$ in $S$. The function $f$ defines a map $ F : U_{2} \to \mathscr{C}^{0}(\bar P, \C)$ where $ \mathscr{C}^{0}(\bar P, \C)$ is the Banach space of continuous functions on $\bar P$, via the formula $F(t)[s] = f(s,t)$ for $t \in U_{2}$ et $s \in \bar P$. The map $F$ is holomorphic: this is an easy consequence of Cauchy's formula on a polydisc $U \subset\subset U_{2} $ with fixed $s \in S$ which computes the partial derivatives in $t := (t_{1}, \dots, t_{n})$:
$$
 \frac{\partial f}{\partial t_{i}}(s, t) = \frac{1}{(2i\pi)^{n}}\int_{\partial\partial U} f(s,\tau).\frac{d\tau_{1}\wedge \dots \wedge d\tau_{n}}{(\tau_{1} - t_{1})\dots (\tau_{i}-t_{i})^{2}\dots (\tau_{n}-t_{n})} \quad \forall t \in U \quad \forall i \in [1,n].
$$
This shows that $F$ is $\C-$differentiable and its differential in $t \in U$ is given by  $h \mapsto \sum_{i=1}^{n} F_{i}(t).h_{i}, \ h \in \C^{n}$, where $F_{i} $ is the map associated to the function
$$ (s,t) \mapsto  \frac{\partial f}{\partial t_{i}}(s, t) \quad i \in [1,n]$$
which is holomorphic for any fixed $s \in S$ thanks to the Cauchy formula above.\\
Let  $H(\bar P,\C)$ be the (closed) subspace of  $\mathscr{C}^{0}(\bar P, \C)$ of continuous functions which are holomorphic on $P$. Our assumption implies that the restriction of $F$ to $U_{1}$ takes its values in this subspace. Let us show that for each point $t \in U_{2}, F(t)$ is still in  $H(\bar P,\C)$: assume this is not true. Then there exists $t_{0} \in U_{2}$ with $F(t_{0})\not\in H(\bar P,\C)$, and so, by  the Hahn-Banach theorem, there exists a continuous linear form  $\lambda$  on $\mathscr{C}^{0}(\bar P, \C)$, vanishing on  $H(\bar P,\C)$ and such that $\lambda(F(t_{0})) \not= 0$. But the  function  $t \mapsto \lambda(F(t))$ is holomorphic on $U_{2}$ and vanishes on $U_{1}$. So it vanishes identically contradicting the fact that $\lambda(F(t_{0})) \not= 0$. So $F$ is an holomorphic map with values in  $H(\bar P,\C)$ and  $f$ is holomorphic on $S\times U_{2}$ when $S$ is smooth.\\
The case where $S$ is a weakly normal complex space is then an immediate consequence, as the continuity of $f$ on $S \times U_{2}$ and the holomorphie of $f$ on $S_{reg}\times U_{2}$, obtained above, are enough to conclude.\\
When $S$ is a general reduced complex space the function $f$ is then a continuous meromorphic function on $S\times U_{2}$ which is holomorphic on $S \times U_{1}$. So the closed analytic subset $Y \subset S \times U_{2}$ of points at which $f$ is not holomorphic has empty interior in each $\{s\}\times U_{2}$. So the criterium  3.1.7 of analytic continuation of chapter IV  in [B-M] allows to conclude.$\hfill \blacksquare$\\

\parag{Remark} It is an easy exercise to weaken the hypothesis of the previous proposition replacing the continuity of $f$ by the hypothesis ``$f$ is measurable and locally bounded  on $S \times U_{2}$''. Then replace the Banach space $\mathscr{C}^{0}(\bar P, \C)$ by the Banach space of bounded measurable functions on $\bar P$ and  in the second step consider the case where $S$ is normal; conclude following the same lines.$\hfill \square$

\parag{Proof of the theorem \ref{cycles}} Consider the graph $\vert G\vert \subset S \times M$ of the f-continuous family  $(X_{s})_{s\in S}$ and let $A$ be the open set of points in $\vert G\vert$ such that the family is analytic in a neighborhood. Precisely, the point $(\sigma,\zeta) \in \vert G\vert$ is in $A$ if there exist open neighborhoods $S_{\sigma}$ and $U_{\zeta}$ respectively of $\sigma$ in $S$ and of  $\zeta$ in $M$ such that the family $(X_{s} \cap M_{\zeta})_{s \in S_{\sigma}}$ is analytic. Remark that, because of our assumption,  $A$ meets each irreducible component of  $\{s_{0}\}\times \vert X_{s_{0}}\vert$.\\
First, assume that there exists a smooth point of $\vert X_{s_{0}}\vert$  in the boundary of the set  $A \cap (\{s_{0}\}\times \vert X_{s_{0}}\vert)$.
Consider now such a point $(s_{0}, z_{0})$ and choose also  a $n-$scale $E := (U,B,j)$ adapted to $X_{s_{0}}$ satisfying the following conditions:
 \begin{align*}
 & \deg_{E}(\vert X_{s_{0}}\vert) = 1 \quad
 & j_{*}(X_{s_{0}}) = k.(U \times \{0\})\quad
  & z_{0} \in j^{-1}(U\times B) \quad
  & j(z_{0}) := (t_{0}, 0).
 \end{align*}
 Then we have a continuous classifying map $f : S_{1}\times U \to \Sym^{k}(B)$ where $S_{1}$ is an open neighborhood of $s_{0}$ in $S$. The map $f$ is holomorphic for each fixed $s \in S_{1}$. As the point  $(s_{0}, z_{0})$ is in the boundary of the open set $A \cap (\{s_{0}\}\times \vert X_{s_{0}}\vert)$ of $\{s_{0}\}\times \vert X_{s_{0}}\vert$, there exists a (non empty) polydisc $U_{1} \subset U$ such that the restriction of $f$ to $S_{1}\times U_{1}$ is holomorphic near $s_{0}$. So, up to shrink $S_{1}$, we can assume that $f$ is holomorphic on $S_{1}\times U_{1}$. Applying the proposition \ref{fonction} to each scalar component of $f$, we conclude that $f$ is holomorphic on $S_{1}\times U$. As we can apply the previous argument to any linear projection of $U\times B$ on $U$ near the vertical one, we obtain also the fact that $f$ is an isotropic map. This contradicts the fact that $(s_{0}, z_{0})$ is in the boundary of $A \cap (\{s_{0}\}\times \vert X_{s_{0}}\vert)$.\\
 If the boundary of $A \cap (\{s_{0}\}\times \vert X_{s_{0}}\vert)$ is contained in the singular set of $\vert X_{s_{0}}\vert$, then we can apply the criterium ([B-M]  ch.IV  crit\`{e}re  3.9.1)  to obtain directly that $A$ contains $\vert X_{s_{0}}\vert$.\\
 So in all cases, the family $(X_{s})_{s\in S}$ is analytic at $s_{0}$. As the graph $\vert G\vert$ is quasi-proper on $S$ by assumption, it is enough to apply the proposition \ref{open} to conclude the proof.$\hfill \blacksquare$\\
 
  \begin{prop}\label{open}
 Let $M$ and  $S$ be a reduced complex spaces and let $(X_{s})_{s \in S}$ be a f-continuous family of $n-$cycles in $M$. Assume that this family is analytic at $s_{0}$. Then there exists an open  neighbourhood $S'$ of $s_{0}$ in $S$ such the family $(X_{s})_{s \in S'}$ is a f-analytic family of $n-$cycles in $M$.
 \end{prop}

 \parag{Proof} The quasi-properness hypothesis on the graph of the family gives an open neighbourhood $S_{1}$ of $s_{0}$ in $S$ and a relatively compact open set $W$ in $M$ such that any irreducible component of any $X_{s}, s \in S_{1}$, meets $W$.  Now for each irreducible component $\Gamma$ of $\vert X_{s_{0}}\vert$ choose a point $z_{\Gamma}\in \Gamma \cap W$, an open neighbourhood $S_{\Gamma}$ of $s_{0}$ in $S_{1}$ and an open neighbourhood $W_{\Gamma}$ of $z_{\Gamma}$ in $W$ such that the family $(X_{s}\cap W_{\Gamma})_{s \in S_{\Gamma}}$ is analytic. Then, as there is a finite number of $\Gamma$, the subset $S' := \cap_{\Gamma} S_{\Gamma}$ is an open neighbourhood of $s_{0}$ in $S_{1}$, and the family $(X_{s}\cap V)_{s \in S'}$ is analytic where $V := \cup_{\Gamma} \ W_{\Gamma} $. Then, for each point $s_{1}\in S'$, we may apply the first part of the proof of the  theorem \ref{cycles}  which is given above (without the present proposition) and conclude that our family is holomorphic at  any point $s_{1}\in S'$. $\hfill \blacksquare$\\
 
 Remark that the previous result is not true in general for a  family of cycles which is analytic at $s_{0}$  and is not $f-$continuous in a neighbourhood of $s_{0}$.

  \subsection{The semi-proper direct image theorem.}

The aim of this paragraph is to give an improvement to the  following  important result given by the theorem 5.0.5 of [B.13]. First recall this result from {\it loc. cit.}
  
 \begin{thm}\label{semi-proper direct image}
 Let $M$  and $S$ be a reduced complex spaces and fix an integer $n$. Assume that we have a holomorphic map $\varphi : S \to \mathcal{C}_{n}^{f}(M)$ which is semi-proper\footnote{This means that for  any cycle $C \in \mathcal{C}_{n}^{f}(M)$ we can find an open neighborhood $\mathcal{V}$ of $C$ and a compact set $K$ in $S$ such that
  \ $ \varphi(S) \cap \mathcal{V} = \varphi(K) \cap \mathcal{V}$.}. Then the closed subset $\varphi(S)$ has a natural  structure of weakly normal complex space such that the tautological family of $n-$cycles parametrized by $\varphi(S)$ is an f-analytic family of cycles in $M$. 
    \end{thm}
    
   Here is our improvement.
    
     \begin{thm}\label{semi-proper direct image bis}
 Let $M$  and $S$ be a reduced complex spaces and fix an integer $n$. Assume that we have a holomorphic map $\varphi : S \to \mathcal{C}_{n}^{f}(M)$ which is semi-proper. Then the closed subset $\varphi(S)$ in $ \mathcal{C}_{n}^{f}(M)$ is analytic and  the sheaf  induced by the sheaf of holomorphic functions on $ \mathcal{C}_{n}^{f}(M)$ defines a structure of reduced (locally finite dimensionnal) complex space on it. 
 \end{thm}
  
 This improvement of the theorem 5.0.5 of [B.13] will allow us to avoid to restrict several statements to  weakly normal complex spaces. It is also an opportunity to precise the proof of this delicate result.\\
 
 The main tool for the proof of the theorem \ref{semi-proper direct image bis} is the semi-proper direct image theorem of [M.00] and the theorem \ref{cycles}.

Before going to the proof of the theorem \ref{semi-proper direct image bis} we have to give some details here on a crucial point for the use of the semi-proper direct image theorem with value in a Banach space in order to have a semi-proper direct image theorem with values in $\mathcal{C}_{n}^{f}(M)$ for any complex space $M$. This is to be compared with condition $(H)$ in [M.00].

\begin{prop}\label{H}
Let $M$ and $S$ be reduced complex spaces and $\varphi : S \to \mathcal{C}_{n}^{f}(M)$ be a holomorphic semi-proper map. Then for each $C \in \varphi(S)$ there exists an open neighborhood $\mathcal{V}$ of $C$ in $\mathcal{C}_{n}^{f}(M)$ and a holomorphic map $f : \mathcal{V} \to \mathcal{U}$ in an open set $\mathcal{U}$ of a Banach space, such that the map $ f\circ\varphi : \varphi^{-1}(\mathcal{V}) \to \mathcal{U}$ is semi-proper.
\end{prop}

\parag{Proof} First recall that a basis of the topology of $\mathcal{C}_{n}^{f}(M)$ is given by finite intersections of  open sets of the type  $\Omega_{k}(E)$ and $\Omega(W)$; here   $\Omega_{k}(E)$, associated to an natural integer $k \geq 0$ and a $n-$scale $E$ on $M$, is defined as the subset of $C \in  \mathcal{C}_{n}^{f}(M)$ such that $E$ is adapted to $C$ and $\deg_{E}(C) = k$ and the open set $\Omega(W)$, associated to a relatively compact open set $W$ in $M$,  is defined as the subset of $C \in  \mathcal{C}_{n}^{f}(M)$ such that  each irreducible component of $C$ meets $W$.\\

 Fix $C_{0} \in \mathcal{C}_{n}^{f}(M)$ and, using the semi-properness of $\varphi$,  let 
 $$\mathcal{V} := \big(\cap_{j=1}^{J} \Omega(W_{j})\big)\cap \big(\cap_{i=1}^{N} \Omega_{k_{i}}(E_{i})\big)$$
  an open neighborhood of $C_{0}$ in $ \mathcal{C}_{n}^{f}(M)$ and $K \subset S$ a compact set in $S$ such that
  \begin{equation*}
   \varphi(S) \cap \mathcal{V} = \varphi(K)\cap \mathcal{V}. \tag{*}
   \end{equation*}
  As $\varphi(K)$ is a compact subset in $ \mathcal{C}_{n}^{f}(M)$, there exists a relatively compact open set $W_{0}$ in $M$ such that for any $s \in K$ any irreducible component of $\varphi(s)$ meets $W_{0}$. As we can make $\mathcal{V}$ smaller around $C_{0}$ keeping the condition $(^{*})$ above, we shall assume that $W_{0}$ is one of the $W_{i}, i \in [1,M]$ and that the union of the centers of the scales $E_{i}, i \in [1,N]$ contains the compact set $\bar W_{0}$.\\
  We consider now the holomorphic map
  $$ f : \mathcal{V} \cap \varphi(S) \longrightarrow \prod_{i=1}^{N} H(\bar U_{i}, \Sym^{k_{i}}(B_{i})) \subset \mathcal{U} $$
  given by $f(C) = (f_{i})_{ i \in [1,N]}$ where $f_{i} \in H(\bar U_{i}, \Sym^{k_{i}}(B_{i}))$ is the classifying map of $C$ in the scale $E_{i}:=(U_{i}, B_{i}, j_{i})$ and where $\mathcal{U}$ is an open set in the Banach space $\prod_{i=1}^{N} H(\bar U_{i}, E_{i}(k_{i}))$\footnote{We denote by $E_{i}(k_{i})$  the Zariski tangent space at $k.\{0\}$ of $\Sym^{k_{i}}(\C^{p_{i}})$. So $\Sym^{k_{i}}(\C^{p_{i}})$ is an algebraic subset of $E_{i}(k_{i})$ ; see [B-M]  ch.I.} such that $\prod_{i=1}^{N}  H(\bar U_{i}, \Sym^{k_{i}}(B_{i}))$ is a closed analytic subset of $\mathcal{U}$. Then $f : \mathcal{V} \to \mathcal{U}$ is holomorphic and injective. We shall show that $f$ induces an homeomorphism of $\mathcal{V}\cap\varphi(S)$ onto $f(\mathcal{V}\cap\varphi(S))$.\\
   As this map is bijective and continuous, it is enough to prove that this map is closed; if  the sequence $(f(C_{\nu}))$ converges to $f(C)$ in $\mathcal{U}$ where $C_{\nu}$ and $C$ are in $\mathcal{V}\cap \varphi(S)$, we shall prove that the sequence $(C_{\nu})$ converges to $C$ in $\mathcal{V}$ (that is to say in $ \mathcal{C}_{n}^{f}(M)$).\\
  The convergence of $(f(C_{\nu}))$ to $f(C)$ implies the convergence of  $(C_{\nu}\cap W_{0})$ to $C \cap W_{0}$ in $\mathcal{C}_{n}^{loc}(W_{0})$,  as we assume that the union of the centers of the scales $E_{i}$ covers $\bar W_{0}$ (see [B-M] ch. IV). Write $C_{\nu} = \varphi(s_{\nu})$ where $s_{\nu}$ is in $K$. Up to pass to a subsequence, we can assume that the sequence  $(s_{\nu})$ converges to $t \in K$. Then $(C_{\nu})$ converges to $D := \varphi(t)$ in $\mathcal{C}_{n}^{f}(M)$ as $\varphi$ is continuous. Then $D$ is in $\varphi(K)$ but it is not, a priori, clear that $D$ is in $\mathcal{V}$ and then equal to $C$. Nevertheless we have $C \cap W_{0} = D \cap W_{0}$ by the uniqueness of a limit in $\mathcal{C}_{n}^{loc}(W_{0})$. As $C$ and $D$ are in $\varphi(K)$, each irreducible components of $C$ and $D$ meets $W_{0}$, and this allows to conclude that $D = C$. As each converging subsequence of the sequence $(C_{\nu})$ converges to $C$ and as we are in the compact set $\varphi(K)$ of $\mathcal{C}_{n}^{f}(M)$, we obtain that the sequence $(C_{\nu})$ converges to $C$ in $\mathcal{C}_{n}^{f}(M)$.\\
  Now the composed map $f\circ\varphi : \varphi^{-1}(\mathcal{V}) \to f(\mathcal{V}\cap\varphi(S)) $ is semi-proper. Then, up to shrink the open set  $\mathcal{U}$ around $f(\mathcal{V}\cap\varphi(S)) $, the map $f\circ\varphi : \varphi^{-1}(\mathcal{V}) \to \mathcal{U}$ is also semi-proper\footnote{The image of a semi-proper map is locally compact;  so it is locally closed in $\mathcal{U}$.} .

  \parag{Proof of the theorem \ref{semi-proper direct image bis}} Assume now that in the previous proof we choose for each  scale $E_{i}, i \in [1,N]$ an open polydisc $U'_{i} \subset\subset U_{i}$ such that the union of the $j_{i}^{-1}(U'_{i}\times B_{i})$ covers the compact set $\bar W_{0}$ and that we replace the Banach analytic set $H(\bar U_{i}, \Sym^{k}(B_{i}))$ by the Banach analytic set  $\Sigma_{U_{i}, U'_{i}}(k_{i})$ which classifies the $f \in H(\bar U_{i}, \Sym^{k}(B_{i}))$ which are isotropic on $\bar U'_{i}$. As the natural projection $ \Sigma_{U_{i}, U'_{i}}(k_{i}) \to H(\bar U_{i}, \Sym^{k}(B_{i}))$ is an holomorphic homeomorphism (see [B.75] chapter III or  [B-M 2] chapter V), and as the map
  $$ \tilde{f} : \mathcal{V} \to \prod_{i=1}^{N} \Sigma_{U_{i}, U'_{i}}(k_{i}) $$
 is holomorphic, we obtain that the map $\tilde{f}\circ\varphi : \varphi^{-1}(\mathcal{V}) \to  \prod_{i=1}^{N} \Sigma_{U_{i}, U'_{i}}(k_{i})  $ is semi-proper on its image and so there exists an open set $\tilde{\mathcal{U}}$ in the ambient Banach space of $ \prod_{i=1}^{N} \Sigma_{U_{i}, U'_{i}}(k_{i}) $ such the map $ \tilde{f} : \mathcal{V} \to \tilde{\mathcal{U}}$ is semi-proper. So the semi-proper direct image theorem with values in a Banach open set gives that $ \tilde{f}(\mathcal{V} \cap \varphi(S))$ is a locally  finite dimensional analytic subset in $\tilde{\mathcal{U}}$.\\
 Consider now the f-continuous  family of $n-$cycles in $M$  parametrized by $ \tilde{f}(\mathcal{V} \cap \varphi(S))$ induced via the homeomorphism $\tilde{f}$ by the tautological family on $\mathcal{V} \cap \varphi(S)$. The universal properties of the Banach analytic spaces  \ $ \Sigma_{U_{i}, U'_{i}}(k_{i})$ imply that this family of $n-$cycles in $M$ is analytic on the open set $\cup_{i=1}^{N} j_{i}^{-1}(U'_{i}\times B_{i})$ which contains the open set $W_{0}$. As $W_{0}$ meets each irreducible component of each cycle in $\mathcal{V} \cap \varphi(S)$, we can apply the theorem \ref{cycles} and conclude that this family is f-analytic. By the universal property\footnote{Recall that we have already proved that $ \tilde{f}(\mathcal{V} \cap \varphi(S))$ is a reduced (locally finite dimensional) complex space.} of $\mathcal{C}_{n}^{f}(M)$, the corresponding classifying map
  $$ \tilde{f}^{-1} :  \tilde{f}(\mathcal{V} \cap \varphi(S)) \longrightarrow  \mathcal{V} \cap \varphi(S) \subset \mathcal{C}_{n}^{f}(M)$$
  is holomorphic. To conclude, it is enough to use the fact that  $ \tilde{f}$ and $ \tilde{f}^{-1}$ are holomorphic between $\mathcal{V} \cap \varphi(S)$ and $ \tilde{f}(\mathcal{V} \cap \varphi(S))$.$\hfill \blacksquare$\\ 
  
  \parag{Remark} We can avoid to use the full strength of the theorem \ref{cycles} in the proof above, because the semi-properness of the map $\tilde{f}\circ \varphi : \varphi^{-1}(\mathcal{V}) \to  \tilde{f}(\mathcal{V} \cap \varphi(S))$ implies the semi-properness of the holomorphic map
  $$(\tilde{f}\circ \varphi) \times \id_{M} : \vert G\vert_{\vert \varphi^{-1}(\mathcal{V})} \to  \tilde{f}(\mathcal{V} \cap \varphi(S))\times M $$ 
  where $\vert G\vert \subset S \times M$ is the set-theoretic graph of the f-analytic family classified by $\varphi$. Then the (finite dimensional) semi-proper direct image theorem of Kuhlmann gives the analyticity of the set-theoretic graph of the f-continuous family parametrized by $\tilde{f}(\mathcal{V} \cap \varphi(S))$, so the fact that this family is weakly holomorphic.$\hfill \square$
 
 \subsection{Meromorphic maps to $\mathcal{C}_{n}^{f}(M)$.}
 
 \parag{Terminology}\begin{itemize}
 
 \item  A {\bf modification} between two reduced complex spaces will be always a proper holomorphic map which induces an isomorphism between two dense Zariski open sets.
  \item We shall say that a holomorphic map  $p : P \to N$ between two irreducible complex spaces is {\bf dominant} if there exists an open dense subset $\Lambda \subset P$ such that  the restriction of $p$ to $\Lambda$ is an open map with dense image in $N$. When $P$ is not irreducible, we shall say that $p$ is dominant when its restriction to each irreducible component of $P$ is dominant.
\end{itemize}

Remark that this implies that the restriction of $p$ to the smooth parts of $P$  has generically maximal rank equal to $\dim N$.

 \begin{defn}\label{strict fiber product}
 Let $\pi : M \to N$ and $p : P \to N$ be two holomorphic  maps where $M, N, P$ are  complex spaces. Assume that $N$ is irreducible and  that $\pi$  and $p$ are dominant. Then we shall call the {\bf strict fiber product} of $\pi$ and $p$,  the map
 $$ \tilde{\pi} :   \tilde{M} := M\times_{N,str}P \to P$$
 which is the restriction of the projection of the usual fiber product to the union of the irreducible components of $M\times_{N} P$ which are dominant  over $P$.
 \end{defn}
 
  Of course, the strict fiber product has two holomorphic projections on $M$ and $P$ which factorize its  natural projection on $N$. They are both dominant. \\
 
 The reader will find easily examples of fiber products which are not equal to the corresponding strict fiber product. For instance, if $\tau : M \to N$ is a modification which is not injective, the fiber product $M\times_{N} M$ has at least one  irreducible component which is not contained in the corresponding  strict fiber product which coincides with the map $\tau$ itself.\\
 
It is an easy exercise left to the reader to prove that if we assume that $p : P \to N$ is a modification and that $\pi : M \to N$ is dominant  the projection $M \times_{N,str}P \to M$ is a modification.\\
  
 The definition of a meromorphic map from a reduced complex space $N$ to $\mathcal{C}_{n}^{f}(M)$ is the usual one.
 
 \begin{defn}\label{mero}
 Fix a complex space $M$ and an integer $n$. Let $N$ be a reduced complex space and $\Sigma \subset N$ a nowhere dense closed analytic subset in $N$. We shall say that a holomorphic map $\varphi : N \setminus \Sigma \to \mathcal{C}_{n}^{f}(M)$ is {\bf meromorphic along $\Sigma$} (or more simply that $\varphi : N ---> \mathcal{C}_{n}^{f}(M)$ is meromorphic) when there exists a modification $\sigma : N_{1} \to N$ with center in $\Sigma$ and a holomorphic map $\varphi_{1} : N_{1} \to \mathcal{C}_{n}^{f}(M)$ extending $\varphi\circ \sigma$.
\end{defn}

To define the ``graph'' of such a meromorphic map we shall need the following corollary of the theorem \ref{semi-proper direct image bis}.

\begin{cor}\label{appli.}
 Fix a complex space $M$ and an integer $n$. Let $N$ and $P$ be reduced complex spaces and let $\varphi : N \to P \times \mathcal{C}_{n}^{f}(M)$ be a semi-proper holomorphic map. Then $\varphi(N)$ is a closed analytic subset in $P \times \mathcal{C}_{n}^{f}(M)$ which is locally of finite dimension.
 \end{cor}
 
 The proof will be an easy consequence of the theorem \ref{semi-proper direct image bis} using the following lemma.
 
 \begin{lemma}
  Fix a complex space $M$ and an integer $n$. Let $P$ be a reduced complex space. Denote $p : P\times M \to P$ and $q : P\times M \to M$ the projections. Then the closed analytic subset $\mathcal{C}_{n}^{f}(p) \subset \mathcal{C}_{n}^{f}(P \times M)$ is bi-holomorphic to the product $P \times \mathcal{C}_{n}^{f}(M)$.
  \end{lemma}
  
  \parag{Proof} Denote $\alpha : \mathcal{C}_{n}^{f}(p) \to P$ the natural holomorphic projection (see the proposition \ref{relatif}) and $\beta : \mathcal{C}_{n}^{f}(p) \to \mathcal{C}_{n}^{f}(M)$ be the holomorphic map induced by the direct image by $q$. Remark that for any $p-$relative $n-$cycle the restriction of $q$ is a holomorphic homeomorphism on its image which is closed (so the restriction of $q$ to such a cycle  is a homeomorphism), then the direct image theorem with parameter applies (see [B.75] or  [B-M] ch.IV). So the map $(\alpha, \beta) : \mathcal{C}_{n}^{f}(p)  \to P \times \mathcal{C}_{n}^{f}(M)$ is holomorphic and bijective. The inverse map $\gamma$,  given by $\gamma(p, C) := \{p\}\times C \in \mathcal{C}_{n}^{f}(p) \subset \mathcal{C}_{n}^{f}(P \times M)$, is also holomorphic  thanks to the product theorem\footnote{Note that even in this case the product theorem is not obvious because a $n-$scale on $P \times M$ adapted to a cycle like  $\{p\}\times C$  does not necessary comes from a $n-$scale on $M$.} for analytic families of cycles (see also [B.75] or [B-M] ch.IV).$\hfill \blacksquare$
  
  \parag{Proof of the corollary} The holomorphic map $\varphi$ defines a holomorphic map $\psi : N \to \mathcal{C}_{n}^{f}(P \times M)$ with values in $ \mathcal{C}_{n}^{f}(p)$. The map $\psi$ is holomorphic and semi-proper and we conclude by applying the theorem \ref{semi-proper direct image bis}.$\hfill \blacksquare$
  
  \begin{lemma}\label{graph mero 1}
  Let $\varphi :  N ---> \mathcal{C}_{n}^{f}(M)$ be a meromorphic map and consider a modification $\sigma : N_{1} \to N$   such that there exists a holomorphic map $\varphi_{1} : N_{1} \to \mathcal{C}_{n}^{f}(M)$ extending $\varphi\circ \sigma$. Then the holomorphic map $(\sigma, \varphi_{1}) : N_{1} \to N \times \mathcal{C}_{n}^{f}(M)$ is proper and its image is the closure in $N \times \mathcal{C}_{n}^{f}(M)$ of the graph of $\varphi_{\vert N \setminus \Sigma}$.
  \end{lemma}
  \parag{Proof} To prove the properness of the map $(\sigma, \varphi_{1})$ we have to prove that it is a closed map and that all fibers are compact. The second point is obvious as $\sigma$ is proper. So consider a closed set $F$ in $N$ and a sequence $(y_{\nu})$ in $F$ such that $(\sigma(y_{\nu}), \varphi_{1}(y_{\nu}))$ converges to $(x, C) \in N \times \mathcal{C}_{n}^{f}(M)$. Take a compact neighbourhood $V$  of $x$ in $N$. As $\sigma$ is proper, $\sigma^{-1}(V)$ is compact and contains all $y_{\nu}$ for $\nu$ large enough. So, up to pass to a subsequence, we can assume that $(y_{\nu})$ converges to $y$ in $N_{1}$. Then the sequence $(\varphi_{1}(y_{\nu}))$ converges to $\varphi_{1}(y)$ in $\mathcal{C}_{n}^{f}(M)$. So we have $y \in F$ and $(x, C) = (\sigma, \varphi_{1})(y)$. The last assertion is obvious.$\hfill \blacksquare$\\
  
  Note that the projection $\tau$  of the reduced complex space $\tilde{N} := (\sigma, \varphi_{1})(N_{1})$ on $N$ is a proper modification and that the map $\tilde{\varphi} : \tilde{N} \to \mathcal{C}_{n}^{f}(M)$ is holomorphic and extends $\varphi\circ \tau$. Moreover, for any modification $\sigma_{2} : N_{2} \to N$ such that there exists a holomorphic map $\varphi_{2} : N_{2} \to \mathcal{C}_{n}^{f}(M)$ extending $\varphi\circ \sigma_{2}$ the holomorphic map $(\sigma_{2}, \varphi_{2})$ factorizes through $\tilde{N}$, meaning that there exists a holomorphic map $h : N_{2} \to \tilde{N}$ such that $\sigma_{2} = \tau\circ h$ and $\varphi_{2} = \tilde{\varphi}\circ h$.

   \begin{defn}\label{graph mero 2}
   In the situation of the lemma above the reduced complex space $\tilde{N} \subset N \times  \mathcal{C}_{n}^{f}(M)$ will be called the {\bf graph} of the meromorphic map $\varphi$.
   \end{defn}

\section{Strongly quasi-proper maps.} 

We first explain that geometrically f-flat maps (resp. strongly quasi-proper maps) is the class of holomorphic maps admitting a holomorphic fiber map (resp. a meromorphic fiber map) defined on the target space  with values in the space of finite type cycles  of the source space (of the appropriate dimension). Then we give some basic results for these notions including stability results in order to dispose of easy criteria ensuring that a holomorphic map is strongly quasi-proper, condition which is not so simple to verify directly on the definition.\\
Of course these tools are essential for the applications given in [B.13] mainly because they allow to use the holomorphic semi-proper direct image theorem \ref{semi-proper direct image bis} with values in the space $\mathcal{C}_{n}^{f}(M)$ where $M$ is any complex space and $n$ any integer. for an example of  how these tools allow to give some results on meromorphic quotients for some holomorphic action of a complex Lie group on a reduced complex space see [B.15]. \\
    
\subsection{Geometrically f-flat maps.}

We recall first the notion of geometrically f-flat  holomorphic map.

\begin{defn}\label{GF} Let $M$ be a pure dimensional complex space and let  $N$ be an irreducible complex space. Put $n := \dim M - \dim N$. Let $\pi : M \to N$ be a holomorphic  map. We shall say that $\pi$ is a {\bf geometrically f-flat} map ({\bf GF map} for short) if there exists a holomorphic map $\varphi : N \to \mathcal{C}_{n}^{f}(M)$ such that for $y$ generic\footnote{One $y$ such $\varphi(y)$ is reduced is enough as $N$ is irreducible ; see lemma \ref{reduced}.} in $N$ the cycle $\varphi(y)$ is reduced and such that the projection on $M$ of the graph $G$ of the f-analytic family of $n-$cycles in $M$  classified  by $\varphi$, is  an isomorphism. 
$$\xymatrix{M \ar[rd]^{\pi} \ar[rr]^{\simeq}& & G \ar[ld]_{pr} \subset N \times M \\ & N & }$$
\end{defn}

Note that the lemma \ref{finitude comp.} implies that in this situation $M$ has finitely many irreducible components. The restriction of $\pi$ to each of these irreducible components is quasi-proper and equidimensional (so surjective), but it is not true that these restrictions are GF maps ; see the example below. But it is not far to be true as we shall see thanks to the  lemma \ref{wGF 1} and the remark 2 following the definition \ref{wGF}

\parag{Example} Let $X := \{ (x,y) \in \C^{2} \ / \  x^{2} = y^{3}\}$ and let $f : X \to \C$ the continuous meromorphic function defined by $f(x,y) = x/y $. Let $G_{\pm} \subset X \times \C$ the graph of $\pm f$ and define $Y := G_{+} \cup G_{-} \subset X \times \C$. Then the projection $\pi : Y \to X$ is a GF map, as the map $F : X \to \Sym^{2}(\C) \simeq \C^{2}$ given by $F(x,y) = (0, y)$ is holomorphic and classifies the fibers of $\pi$ (note that $f^{2}(x,y) = y$). But the restriction of the projection to the irreducible component $G_{+}$ of  \ $Y$ \ is not a GF map as $f$ is not holomorphic on $X$.\\

\begin{lemma}\label{fiber GF}
Let $\pi : M \to N$ be a holomorphic  $n-$equidimensional map between a pure dimensional complex space $M$ to an irreducible complex space $N$, where the integer $n$  is equal to $ \dim M - \dim N$.  Assume that there exists a holomorphic map 
 $$\varphi : N \to \mathcal{C}_{n}^{f}(M)$$
 such that for any $y$ in a  dense set in $N$ the cycle $\varphi(y)$ is reduced and equal to the fiber of $\pi$. Then $\pi$ is a GF map.
\end{lemma}

\parag{Proof} As for $y$ in a dense set we have $\varphi(y) \in \mathcal{C}_{n}^{f}(\pi)$ which is a closed (analytic) subset in $\mathcal{C}_{n}^{f}(M)$ (see the proposition \ref{relatif}), we have $\vert \varphi(y)\vert \subset \pi^{-1}(y)$ for any $y \in N$. Consider the second projection $pr : G \to M$ of the graph $G \subset N \times M$ of the f-analytic family of $n-$cycles classified by $\varphi$. As $\pi$ is an open map, there exists a dense subset of $x$  in $M$ such that $(\pi(x), x)$ is in $G$. So the holomorphic map $(\pi, \id_{M}) : M \to N \times M$ takes its values in $G$ and it gives a holomorphic  inverse to $pr$. So $\pi$ is a GF map.$\hfill \blacksquare$

\parag{Remark}
 A GF-map is quasi-proper, equidimensional and surjective. But a holomorphic quasi-proper equidimensionnal (and then surjective)  map between irreducible complex spaces is not always a GF map : for instance the (weak) normalization of a reduced complex space $N$ is not a GF-map if $N$ is not (weakly) normal. Nevertheless the next lemma explains that  such a map it is not so far to be a GF map. 

\begin{lemma}\label{wGF 1}
Let $\pi : M \to N$ be a holomorphic map between a pure dimensional complex space $M$ and an irreducible complex space $N$. Let $n := \dim M - \dim N$ and let  $\nu : \tilde{N} \to N$ be the normalization of $N$. Then the following properties are equivalent :
\begin{enumerate}[i)]
\item $\pi$ is quasi-proper and $n-$equidimensional.
\item The projection\footnote{ See the definition \ref{strict fiber product}.} $\tilde{\pi} : \tilde{M} := M\times_{N, str}\tilde{N} \to \tilde{N} $ is a GF-map.
\end{enumerate}
\end{lemma}

\parag{Proof} Assume i) ; first let us prove that $\tilde{\pi}$ is quasi-proper. Let $\tilde{K}$ be a compact set in $\tilde{N}$ ; as $K := \nu(\tilde{K})$ is a compact set in $N$ there exists a compact set $L$ in $M$ such that each irreducible component of a $\pi^{-1}(y)$ for $y \in K$ meets $L$. Then the compact set $(L\times \nu^{-1}(K)) \cap (M\times_{N}\tilde{N})$ meets any irreducible component of a $\tilde{\pi}^{-1}(\tilde{y})$ for each $\tilde{y} \in \nu^{-1}(K)$. So $\tilde{\pi}$ is quasi-proper as any irreducible component of a fiber of $\tilde{\pi}$ is an irreducible component of a fiber of $\pi$. The $n-$equidimensionality of $\tilde{\pi}$ is obvious. Now $\tilde{\pi} : M \times_{N,str}\tilde{N} \to \tilde{N}$  is quasi-proper and equidimensional on a normal complex space, so there exists a  holomorphic map (see [B.M] ch.IV) $\varphi : \tilde{N} \to \mathcal{C}_{n}^{f}(\tilde{M})$ such that $\vert \varphi(\tilde{y}) \vert = \tilde{\pi}^{-1}(s)$ for $\tilde{y}$ general in $\tilde{N}$, where we put $\tilde{M} := M\times_{N, str}\tilde{N}$. Using the analyticity  of the direct image of cycles by $\nu$ the lemma \ref{fiber GF} gives that $\tilde{\pi}$ is a GF map.\\
Conversely, assume ii). First we prove the $n-$equidimensionality of $\pi$. Assume there exists a non normal point $y \in N$ such that the fiber $\pi^{-1}(y)$ has an irreducible component $\Gamma$  of dimension strictly bigger than $n$. This irreducible component does not appears in a fiber of $\tilde{\pi}$. We can find a sequence of points in $\pi^{-1}(N_{norm})$ converging to the generic point of $\Gamma$, where $N_{norm}$ is the open dense set of normal points in $N$. As $\tilde{M}$ is a finite modification of $M$ we can assume, up to pass to a sub-sequence,  that this sequence converges in $\tilde{M}$. We then find a point $y$  in $N$ such that the fiber of $\tilde{\pi}$ contains the generic point of $\Gamma$. This contradicts our assumption.\\
Let us  prove  the quasi-properness of $\pi$. Let $K$ be a compact in $N$ ; then $\nu^{-1}(K)$ is a compact in $\tilde{N}$ and so there exists a compact $L$ in $M$ such that any irreducible component of $\tilde{\pi}^{-1}(\tilde{y})$ for $\tilde{y} \in \nu^{-1}(\tilde{K})$ meets $(L\times \nu^{-1}(K))\cap \tilde{M}$. Then any irreducible component of $\pi^{-1}(y)$ for $y \in K$ meets $L$. $\hfill \blacksquare$

\begin{defn}\label{wGF}
A holomorphic map $\pi : M \to N$ between a pure dimensional complex space $M$ to an irreducible complex space $N$ which satisfies the conditions of the lemma above will be called a {\bf weakly geometrically f-flat} map (in short  a {\bf wGF} map).
\end{defn}

\parag{Remarks}\begin{enumerate}
\item As in the situation of the previous lemma  the projection $M\times_{N, str}\tilde{N} \to M$ is a modification, this implies that $M$ has finitely many irreducible components because the lemma \ref{finitude comp.} implies that $\tilde{M} := M \times_{N,str} \tilde{N}$ has only finitely many irreducible components.
\item  If the holomorphic map $\pi : M \to N$ is wGF map  its restriction to any irreducible component of $M$ is again wGF because this restriction is  equidimensional and quasi-proper.
\item Conversely, if $\pi : M \to N$ is holomorphic, if $M$ is pure dimensional with finitely many irreducible components  and if the restriction of $\pi$ on each irreducible component of $M$  is a wGF map, then $\pi$ is a wGF map, again thanks to the lemma \ref{wGF 1}. 
\item If $\pi : M \to N$ is a wGF map then for any closed irreducible analytic subset $P \subset N$ the restriction  $ \pi_{\pi^{-1}(P)} : \pi^{-1}(P) \to P$ is a wGF map.\\
\end{enumerate}

 The first stability results for this notion of wGF map are given by the following lemma.

 \begin{lemma}\label{fiber product}
 Let $ \alpha : M \to N$ be a wGF map and $\sigma : P \to N$  be a dominant  holomorphic map of an irreducible complex space $P$ to $N$. The the natural projection of  the fiber product $\tilde{\alpha} : M\times_{N, str} P \to P$ is a wGF map. If $\alpha$ is a GF map, so is $\tilde{\alpha}$.
 \end{lemma}
 
 \parag{proof}  Assume first that $\alpha$ is geometrically f-flat; its fibers are classified by a holomorphic map
 $$ \varphi : N \to \mathcal{C}_{n}^{f}(M)$$
 where $n := \dim M - \dim N$. Compose this map with $\sigma$; this gives a holomorphic map $\varphi\circ\sigma$ which classifies a f-analytic family of $n-$cycles in $M$ parametrized by $P$. Let $G \subset P \times M$ the graph of this family. As for $y$ generic in $N$ the cycle $\varphi(y)$ is reduced and the map $\sigma$ is dominant, the graph $G $ of the family parametrized by $P$ is reduced and given by
 $$ G := \{(x,z) \in M \times P \ / \  x \in \vert \varphi(\sigma(z))\vert \}.$$
 But, as $\vert \varphi(y)\vert = \alpha^{-1}(y)$ for each $y \in N$, we see that $G = M\times_{N, str} P$ and that for $z$ in a dense set in $P$ the fiber of the projection of  $G$ on $P$ is the fiber of the map $\tilde{\alpha}$. So this map is geometrically f-flat, thanks to the lemma \ref{fiber GF}.\\
 Consider now the case where $\alpha$ is only a wGF  map. So, thanks to the lemma \ref{wGF 1}, if $\nu : \tilde{N} \to N$ is the normalization of $N$, the fiber product $\tilde{M} := M\times_{N, str}\tilde{N}$ has a projection $\tilde{\alpha}$ on $\tilde{N}$ which is a GF map.  Consider now the normalization $\mu : \tilde{P} \to P$. We have the following commutative diagram, where $\tilde{\sigma}$ is a lifting of \ $\sigma\circ\mu$ \ to \ $\tilde{N}$ \ which exists by normality of $\tilde{P}$ as $\sigma$ is dominant.
 $$\xymatrix{\tilde{M} \ar[r]^{\tilde{\alpha}} \ar[d]_{\tilde{\nu}} & \tilde{N}  \ar[d]_{\nu} & \tilde{P} \ar[l]_{\tilde{\sigma}} \ar[d]_{\mu} \\ M \ar[r]^{\alpha} & N & P  \ar[l]_{\sigma}} $$
  Let $\varphi : \tilde{N} \to \mathcal{C}_{n}^{f}(\tilde{M})$ the holomorphic map classifying the fibers of the GF map $\tilde{\alpha}$. Composed with $\tilde{\sigma}$ and with the direct image map $\tilde{\nu}_{*}$ for $n-$cycles (see [B.M] ch. IV case of a proper map) it gives a holomorphic map
  $$ \psi : \tilde{P}\overset{\tilde{\sigma}}{\longrightarrow}\mathcal{C}_{n}^{f}(\tilde{M}) \overset{ \tilde{\nu}_{*}}{\longrightarrow}  \mathcal{C}_{n}^{f}(M)$$
  and it is easy to see that for $z$ in  a general subset in $\tilde{P}$ its value is the reduced cycle $C$ in $M$ where the cycle $C\times\{z\}$ in $M\times_{N, str} \tilde{P}$ which is the fiber of the natural projection $ p : M\times_{N, str} \tilde{P} \to \tilde{P}$, as we have $(M\times_{N, str} P)\times_{P, str} \tilde{P} \simeq M\times_{N, str} \tilde{P}$. So $p$ is a GF map from the first case proved above, and we conclude the proof thanks to the lemma \ref{wGF 1}.$\hfill \blacksquare$
  
  \parag{Remark} Let $\alpha : M \to N$ be a wGF map and  $\sigma : P \to N$ any holomorphic map, it is easy to see that the projection $M \times_{N,str}P \to P$ is again a wGF map  using the lemma \ref{wGF 1} ; this generalizes the previous lemma and the remark 4 following the definition \ref{wGF}. Note that the case of a GF map is not clear because the image of $\sigma$ may be inside the locus of non reduced fibers of $\alpha$. 

\begin{lemma}\label{compose}
Let $\alpha : M \to N$ and $\beta : N \to P$ two wGF maps  between a pure dimensional complex space $M$ and irreducible complex spaces $N$ and $P$. Then $\beta\circ\alpha$ is a wGF map. 
\end{lemma}

\parag{Proof} Put $m := \dim M - \dim N$ and also  $n := \dim N - \dim P$. Fix an irreducible component $\Delta$ of $(\beta\circ\alpha)^{-1}(z)$ for some $z$ in $P$. Then $\alpha$ induces a map $\alpha' : \Delta \to \beta^{-1}(z)$ where $\beta^{-1}(z)$ has pure dimension $n$ and $\alpha'$ has pure dimension $m$ fibers. So $\Delta$ has dimension at most $m+n$, and, as the map  $\beta\circ\alpha$ is surjective between each irreducible component of  $M$ and the irreducible complex space $P$ with $\dim M - \dim P = m+n$, we obtain that $\Delta$ has dimension $m+n$ and also that $\alpha(\Delta)$ is dense in an irreducible component $\Gamma$ of $\beta^{-1}(z)$. So $\beta\circ\alpha$ is equidimensional (and surjective).\\
 We shall show that it is also quasi-proper. Let $z_{0}$ be a point in $P$. The quasi-properness of $\beta$ implies that there exists an open neighborhood $V$ of $z_{0}$ in $P$ and a relatively compact open set $W$ in $N$ such that any irreducible component $\Gamma$ of any fiber $\beta^{-1}(z)$ for $z \in V$ meets $W$. As $W$ is relatively compact and $\alpha$ quasi-proper, there exists a relatively compact open set $U$ in $M$ such any irreducible component $\Gamma$ of a fiber $\alpha^{-1}(y)$ with $y \in W$ meets $U$. Take now $z \in V$ and $\Delta$ an irreducible component of $(\beta\circ\alpha)^{-1}(z)$. Let $\Gamma$ be the irreducible component of $\beta^{-1}(z)$ in which $\alpha(\Delta)$ is dense (see above). As $\Gamma$ meets $W$ the dense set $\alpha(\Delta)$ in $\Gamma$ meets the non empty open set $\Gamma \cap W$ of $\Gamma$. Let $y$ be a point in $\alpha(\Delta) \cap W$ such that there exists a point $x$ in $\Delta$ with $\alpha(x) = y$ and such that $x$ is a smooth point in $(\beta\circ\alpha)^{-1}(z)$. Such a point $x$ exists because we may choose $y$ to be a smooth point in $\beta^{-1}(z)$ and also in the image of the open dense set in $\Delta$ of points which are smooth in $(\beta\circ\alpha)^{-1}(z)$.\\
 Now an irreducible component $\Gamma$ of $\alpha^{-1}(y)$ containing $x$ has to be contained in $\Delta$ as it is contained in $(\beta\circ\alpha)^{-1}(z)$ and contains $x$. But such a $\Gamma$ meets $U$ as $y$ is in $W$. So $\Delta$ meets $U$. This gives the quasi-properness of $\beta\circ\alpha$.$\hfill\blacksquare$\\
 
 The next lemma will be used later on.
 
 \begin{lemma}\label{modif. and GF}
 Let $ \tau : M_{1} \to M$ be a modification and $\alpha : M \to N$ a holomorphic map between irreducible complex spaces. Assume that $\alpha\circ\tau$ is a GF map. Then $\alpha$ is a GF map.
 \end{lemma}
 
 \parag{Proof} Our assumption gives a holomorphic map $ \varphi : N \to \mathcal{C}_{n}(M_{1})$ which classifies the fibers of $\tau\circ\alpha$. Denote by $\Sigma$ the center of $\tau$. Then the corollary \ref{une comp.}  implies that for generic $y \in N$ no irreducible component of  $\vert \varphi(y)\vert$ is contained in $ \tau^{-1}(\Sigma)$ (see the corollary \ref{une comp.}); then the direct image $\tau_{*}(\varphi(y))$ has its support equal to $\alpha^{-1}(y)$. This means that the composition of $\varphi$ with the direct image of $n-$cycle by $\tau$  gives a holomorphic map $\psi : N \to \mathcal{C}_{n}(M)$ such that for generic  $y \in N$ we have $\psi(y) = \vert \psi(y)\vert = \alpha^{-1}(y)$. As $\alpha$ is $n-$equidimensional because $\dim \alpha^{-1}(z) \leq \dim (\alpha\circ\tau)^{-1}(z)$ for each $z \in N$, the lemma \ref{fiber GF} allows to conclude that $\alpha$ is a GF map.$\hfill \blacksquare$\\

 \begin{lemma}\label{finite and GF}
 Let $\alpha : M \to T$ and $p : T \to N$ holomorphic maps between normal complex spaces such that $p$ is proper finite and surjective and $p\circ\alpha$ is GF. Then $\alpha$ is a GF map.
 \end{lemma}
 
 \parag{Proof} As $N$ is normal, it is enough to show that $\alpha$ is quasi-proper and equidimensional. Let $\delta$ be an irreducible component of $\alpha^{-1}(t)$ for some $t \in T$. As $p$ has finite fibers, $\delta$ is an irreducible component of $(p\circ\alpha)^{-1}(p(t))$ and so has pure dimension $n := \dim M - \dim N$. So $\alpha$ is $n-$equidimensional and any irreducible component of the  fiber of $\alpha$ at a point $t \in T$  is an irreducible component of the fiber of $p\circ\alpha$ at the point $p(t)$.\\
  Let $K$ be a compact subset in $T$; then $p(K)$ is compact and there exists a compact set $L$ in $M$ such that for any $y \in p(K)$ any irreducible component $\gamma$ of $(p\circ \alpha)^{-1}(y)$ meets $L$. Now an irreducible component of a fiber at a point $t \in K$ is an irreducible component of a fiber of  $p\circ\alpha$ at the point $p(t) \in p(K)$. So it meets the compact $L$.$\hfill \blacksquare$\\
  
  It is an simple exercice to show an analoguous result for wGF maps instead of GF maps (and the normality assumptions can be omitted).\\

We conclude this paragraph on GF maps by  the following ``embedding theorem''.
  
  \begin{thm}\label{embed}
  Let $\pi : M \to N$ a holomorphic GF map between irreducible complex spaces. Then the   fiber map $\varphi : N \to \mathcal{C}_{n}^{f}(M)$ is a proper holomorphic embedding. 
  \end{thm}
  
  \parag{Proof} The first step of the proof is to show that the map $\varphi$ is proper, that is to say that $\varphi$  is closed with compact fibers. But as this map is clearly injective we only have to prove that for any closed subset $F \subset N$ its image $\varphi(F)$ is closed in $\mathcal{C}_{n}^{f}(M)$. So consider a sequence $(y_{\nu})$ of points in $F$ such that the sequence $(\varphi(y_{\nu}))$ converges to an element $C_{0}\in \mathcal{C}_{n}^{f}(M)$. Let $W \subset\subset M$ a relatively compact open set in $M$ such that each irreducible component of $\vert C_{0} \vert$ meets $W$. As $\Omega(W)$ is an open set in $\mathcal{C}_{n}^{f}(M)$ containing $C_{0}$ the cycle $\varphi(y_{\nu})$ is in $\Omega(W)$ for $\nu$ large enough. So choose for each such $\nu$ a point $x_{\nu}$ in $W \cap \pi^{-1}(y_{\nu})$. Then each $x_{\nu}$ is in the compact set $\bar W \cap \pi^{-1}(F)$. Up to pass to a sub-sequence we may assume that the sequence $(x_{\nu})$ converges to a point $x \in \bar W \cap \pi^{-1}(F)$. Then $y := \pi(x)$ is in $F$ and is the limit of the sequence $(y_{\nu})$, so we have $C_{0} = \varphi(y)$ with $y \in F$. So the map $\varphi$ is closed and then proper.\\
  Now the theorem \ref{semi-proper direct image bis}\footnote{ In the proper case ; for the generalization of Remmert's direct image theorem with value in a open set in a Banach space see [B-M] ch.III th. 7.3.1.} implies that $N_{\pi} := \varphi(N)$ is a closed analytic subset which is an irreducible complex space and $\varphi : N \to N_{\pi}$ is a holomorphic homeomorphism. To complete the proof of the theorem we have to show that any germ of holomorphic function at a point $y \in N$ is the pull-back of a germ of holomorphic function of $N_{\pi}$ at the point $\varphi(y)$.\\
  So fix a point $y \in N$ and a holomorphic germ $f \in \mathcal{O}_{N,y}$. Choose  a smooth point $x$ of $\vert \varphi(y)\vert$. Let $E := (U, B, j)$ be a $n-$scale on $M$ adapted to $\varphi(y)$ such that $x $ is in $j^{-1}(U\times B)$ and such that $\varphi(y) \cap j^{-1}(U\times B) = k. j^{-1}(U \times \{0\})$. Then we can assume that the holomorphic germ $f \in \mathcal{O}_{N,y}$ is define on the open set $\pi( j^{-1}(U\times B))$ (remember that $\pi$ is open as it is a GF map). Now choose a function $\rho \in \mathscr{C}_{c}^{\infty}(U)$ such that $\int_{U} \ \rho(t).dt\wedge d\bar t = 1/k$. Define the  $(n,n)-$form  $\omega :=  f(\pi(j^{-1}(t,x))).\rho(t).dt\wedge d\bar t $ on $U \times B$; it has a $B-$proper support and it is $\bar \partial-$closed. Then it induces by integration on cycles a holomorphic function on the open set $\Omega_{k}(E)$\footnote{see [B.75] ch. IV or [B-M] ch.IV}  of $\mathcal{C}_{n}^{f}(M)$ via the integration map
  $$ \int_{\square} \ \omega : H(\bar U, \Sym^{k}(B)) \to \C \qquad {\rm given \ by} \quad X \mapsto \int_{X} \ \omega.$$
  On any $\varphi(z) = C \in \varphi(N) \cap \Omega_{k}(E)$ this holomorphic function takes the value $f(z)$, and this proves that $\varphi : N \to N_{\pi}$ is an isomorphism of complex spaces.$\hfill \blacksquare$\\

Note that, conversely, if $\varphi : N \to  \mathcal{C}_{n}^{f}(M)$ is a proper holomorphic embedding, a necessary and sufficient condition for $\varphi$ to be the fiber map of a GF homorphic map $\pi : M \to N$ is the fact that for $y$ generic in $N$ the cycle $\varphi(y)$ is reduced and that the projection $p : G \to M$ of the graph $G \subset N \times M$  of the f-analytic family parametrized by $N$ is an isomorphism onto $M$. Of course this implies that for $y \not= y'$ in $N$ the cycles $\varphi(y)$ and $\varphi(y')$ are disjoint.

\subsection{Strongly quasi-proper map.}

Recall now, for the convenience of the reader, the definition of a strongly quasi-proper map  given in [B.13], using the terminology introduced above. Consider the following situation :

\parag{The standard situation} Let $\pi : M \to N$ be a  holomorphic map between a pure dimensional complex space $M$ and an irreducible complex space $N$. Define the integer $n$ as  $n := \dim M - \dim N$. Assume that there is a closed analytic subset $\Sigma \in N$ with no interior point in $N$ such that the restriction $\pi' : M \setminus \pi^{-1}(\Sigma) \to N\setminus \Sigma $ \  is a GF map. Let $\varphi' : N\setminus \Sigma \to \mathcal{C}_{n}^{f}(M)$ the holomorphic map classifying the fibers of $\pi'$ (as $n-$cycles in $M$) and $\Gamma \subset (N\setminus \Sigma)\times M$ the graph of this family.
 \parag{Remark} For a quasi-proper surjective map $\pi : M \to N$ between irreducible complex spaces, with \  $n := \dim M - \dim N$, there always exists a closed analytic subset $\Sigma \subset N$ with no interior point in $N$ such that  the following  restriction of $\pi$  $\pi' : M \setminus \pi^{-1}(\Sigma) \to N\setminus \Sigma $ \  is a GF map (see the end of the proof of the proposition \ref{SPQ dense} for details).

\begin{defn}\label{sqp}
 In the standard situation described above we shall say that $\pi$ is {\bf strongly quasi-proper} (a { \bf SQP map} for short) if and only if the closure $\bar \Gamma$ of  \ $\Gamma$ in $N \times \mathcal{C}_{n}^{f}(M)$ is proper over $N$.
\end{defn}

Of course we shall show that a SQP map is quasi-proper in the usual sense (see the proof of the next proposition).\\

We begin by an improvement of the criterium given in [B.13] in order that a holomorphic map $\pi : M \to N$ between irreducible complex spaces will be SQP.

\begin{prop}\label{SPQ dense}
Let $h : M \to N$ a holomorphic map between a pure dimensional complex space $M$ and an irreducible complex space $N$. Put $n := \dim M - \dim N$. Assume that there exists a  dense subset $\Lambda$ in $N$ such that $h^{-1}(\Lambda)$ is dense in $M$ and such that for each $y \in \Lambda$ the fiber $h^{-1}(y)$ is non empty, reduced\footnote{We mean here reduced as a cycle ; this is equivalent to the fact that the natural structure of complex space on this fiber is generically reduced ``in the algebraic sense''.} and of (pure) dimension $n$ with finitely many irreducible components. Note $\gamma : \Lambda \to \mathcal{C}_{n}^{f}(M)$ the map defined by $\gamma(y) := h^{-1}(y)$ where the $n-$cycle $\gamma(y)$ is reduced. Let $\Gamma$ the graph of the map $\gamma$ and $\bar \Gamma$ the closure of $\Gamma$ in $N\times \mathcal{C}_{n}^{f}(M)$. Our main assumption is now the following:
\begin{itemize}
\item The natural projection $\tau :\bar \Gamma \to N$ is proper.
\end{itemize}
Then the map $h$ is strongly quasi-proper.
\end{prop}

\parag{Proof} The empty $n-$cycle is open (and closed) in $\mathcal{C}_{n}^{f}(M)$ so any sequence converging to it is stationary. As $\Lambda$ is dense in $N$ and $\gamma(y) \not= \emptyset_{n}$ for $y \in \Lambda$, any fiber of $\tau$ cannot be equal to $\{\emptyset_{n}\}$. But $\tau(\bar \Gamma) = N$ as it is closed and contains $\Lambda$. This implies that $h(M) = \tau(\bar \Gamma) = N$ and $h$ is surjective.\\
Our second  step (and the main step in fact) will be the proof  that the map $h$ is quasi-proper. So fix a point $y \in N$ and let $V$ be an open relatively compact neighbourhood of $y$ in $N$. Fix $y' \in V$ and  choose an irreducible component $C$ of $h^{-1}(y')$\footnote{From the surjectivity of $h$ proved above, $h^{-1}(y')$ is not empty.}. Let $x'$ be a (generic) point in $C$ such that $x'$  does not belong to any other irreducible component of $h^{-1}(y')$. Then we can choose a sequence $(x_{\nu})_{\nu \geq 0}$ in $h^{-1}(\Lambda)$ converging to $x'$. For $\nu \gg 1$ we have $h(x_{\nu}) \in V$ so the cycles $\gamma(h(x_{\nu}))$ are in the compact set $p_{2}(\tau^{-1}(\bar V))$ of $\mathcal{C}_{n}^{f}(M)$. Up to pass to a sub-sequence, we can assume that the sequence $(\gamma(h(x_{\nu})))$ converges to a cycle $\delta \in \mathcal{C}_{n}^{f}(M)$. As we have $x_{\nu}\in \gamma(h(x_{\nu}))$ for each $\nu$ we have $x' \in \vert \delta\vert$. By compactness of the subset $p(\tau^{-1}(\bar V))$  of \ $\mathcal{C}_{n}^{f}(M)$, there exists a compact subset $K$ in $M$ such that any irreducible component of any cycle in $p(\tau^{-1}(\bar V))$ meets $K$. So this is the case for each irreducible component of any cycle $\gamma(h(x_{\nu}))$ and of  $\delta$. Let $\delta_{0}$ be an irreducible component of $\delta$ containing $x'$. So $\delta_{0}$ meets $K$. But $\vert \delta\vert$ is contained in $h^{-1}(y')$ because the condition for a $n-$cycle in $\mathcal{C}_{n}^{f}(M)$ to be contained in a fiber of $h$ is a closed condition (see the proposition \ref{relatif}). As $\delta_{0}$ is irreducible, contained in $h^{-1}(y')$ and contains $x'$, this implies $\delta_{0} \subset C$ and so $C$ meets $K$.
So we have proved that for any $y \in N$ there exists an open neighbourhood  $V$ of $y$ in $N$ and a compact set $K$ in $M$ such that for any $y' \in V$ and any irreducible component $C$ of $h^{-1}(y')$ the intersection $C \cap K$ is not empty. This is the definition of the quasi-properness of the map  $h$.\\
Now we conclude that the subset $\Sigma_{0}$ of $N$ of $y $ such that $\dim h^{-1}(y) > n$ is a closed analytic set with no interior point because $\Sigma_{0}$ is the image by $h$  of  the closed analytic subset
$$ Z := \{ x \in M \ / \  \dim_{x}(h^{-1}(h(x))) > n \} $$
 in $M$ which is an union of  irreducible components of  fibers of $h$; so the restriction of $h$ to this analytic subset is still quasi-proper. So Kuhlmann's theorem gives the conclusion. Now let $\Sigma_{1}$ be the closed analytic subset of non normal points in $N$ and let $N' := N \setminus (\Sigma_{0} \cup \Sigma_{1})$. Then the map $h' : h^{-1}(N') \to N'$ is quasi-proper equidimensional on a normal complex space. Then we have a holomorphic map
$$ \varphi : N' \to \mathcal{C}_{n}^{f}(M) $$
classifying the generic fibers of $h'$. Of course, for $y \in \Lambda \cap N'$ we have $\vert \varphi(y)\vert = \gamma(y)$. But the subset  $T$ of  points $y \in N'$ such that $\varphi(y)$ is not a reduced cycle is a closed analytic subset of $N'$ (see the lemma \ref{reduced}) which has no interior point. So the  dense subset $\Lambda$ meets the dense open set $N'\setminus T$  in a dense subset in $N$. This implies that the closure of the graph of $\varphi$ in $N \times \mathcal{C}_{n}^{f}(M) $ is equal to $\bar \Gamma$. Then, by definition,  the map $h$ is strongly quasi-proper. $\hfill \blacksquare$\\

An immediate corollary of this result is a first stability result of SQP maps.

\begin{cor}\label{invariance}
Let $\pi : M \to N$ be a strongly quasi-proper holomorphic map between irreducible complex spaces and $\tau :T \to N$ a holomorphic map from an irreducible complex space $T$ such that $\tau(T)$ is not contained in $A$  the set of points  in $N$  having ''big'' fibers for $\pi$. Then for any irreducible component $\tau^{*}(M)_{i}$  of $(T \times_{N}M)$ which surjects on $T$ the pull-back $\tau^{*}(\pi)_{i} : \tau^{*}(M)_{i} \to T$ is strongly quasi-proper.  $\hfill \blacksquare$\\
\end{cor}

We give now a characterization of SQP maps in term of GF maps, which is variant of the important theorem 2.4.4 of [B.13]; see [M.00] or  the Appendix (section 4) for a complete proof.

\begin{thm}\label{SQP et GF}
Let $\pi : M \to N$ be a surjective holomorphic map between irreducible complex spaces.  The map $\pi$ is  strongly quasi-proper if and only if  there exists a modification $\sigma : \tilde{N} \to N$ such that the strict transform $\tilde{\pi} : \tilde{M} \to \tilde{N}$ by $\sigma$ of the map $\pi$  is a GF map.
\end{thm}

\parag{Proof} Assume first that $\pi$ is a SQP map. Then there exists a modification $\sigma : \tilde{N} \to N$ and a holomorphic map $\tilde{\varphi} : \tilde{N} \to \mathcal{C}_{n}^{f}(M)$ extending holomorphically the map $\varphi : N' \to \mathcal{C}_{n}^{f}(M)$ classifying the generic fibers of $\pi$ (see section 4). Then define the holomorphic map $\psi : \tilde{N} \to  \mathcal{C}_{n}^{f}(\tilde{M})$ by $\psi(\tilde{y}) := \{\tilde{y}\} \times \tilde{\varphi}(y)$ as a cycle in $\tilde{N}\times M$. The holomorphy is consequence of the product theorem for analytic families of cycles (see [B.75] or  [B.M] ch.IV) and these cycles are in $\tilde{M} := M \times_{N, str}\tilde{N}$ because this is true for generic $\tilde{y} \in \tilde{N}$ and so for all $\tilde{y} \in \tilde{N}$ by continuity of $\psi$. Then we can conclude that $\tilde{\pi}$ is a GF map applying the lemma \ref{fiber GF}.\\
Assume now that we have a modification $\sigma : \tilde{N} \to N$ such that the strict transform $\tilde{\pi} : \tilde{M} \to \tilde{N}$ \ of $\pi$ by $\sigma$  is a GF map. Then the composition of the holomorphic map $\psi : \tilde{N} \to \mathcal{C}_{n}^{f}(\tilde{M})$ classifying the fibers of $\tilde{\pi}$ with the direct image of cycles by the modification $\tilde{\sigma} : \tilde{M} \to M$ induced by the projection on $M$ of the fiber product $M\times_{N, str}\tilde{N}$  gives a holomorphic map $\tilde{\varphi} : \tilde{N} \to \mathcal{C}_{n}^{f}(M)$.  It is clear that, if $\Sigma$ is the center of $\sigma$, the restriction of $\tilde{\varphi}$ is the map $\varphi$ on $N\setminus \Sigma$. So, applying the theorem \ref{embed} we see that $\pi$ is a SQP map because the closure of the graph of its fiber map  is  the image of  $\tilde{N}$ which is proper on $N$.$\hfill \blacksquare$\\

Remark that in the previous theorem we can replace GF by wGF thanks to lemma \ref{wGF 1}.

\begin{prop}\label{mero et SQP}
 A quasi-proper surjective holomorphic map between irreducible complex spaces is SQP if and only if its fiber map, defined and holomorphic on a dense Zariski open set in $M$\footnote{It is enough to take the set of normal points in $N$ intersected with the complement of the ``big fibers'' subset as in the  remark following the definition of the standard situation (see the begining of the paragraph 3.2).}, is meromorphic. Moreover, if $\pi : M \to N$ is a quasi-proper surjective map such that its fiber map is meromorphic, then the strict transform $\tilde{\pi} : \tilde{M} \to N_{\pi}$ of $\pi$ by the modification $\tau : N_{\pi} \to N$ given by the graph of the meromorphic fiber map (see the paragraph 1.4) is a GF map.
  \end{prop}
 
 \parag{Proof} Consider a SQP holomorphic map $\pi : M \to N$. Then, thanks to the theorem \ref{SQP et GF}, there exists a modification $\sigma : \tilde{N} \to N$ such that  the strict transform $\tilde{\pi} : \tilde{M} \to \tilde{N}$  of $\pi$ by  $\sigma$  is a GF map. So there exists a holomorphic map $\tilde{\varphi} : \tilde{N} \to \mathcal{C}_{n}^{f}(\tilde{M})$ classifying the fibers of $\tilde{\pi}$.  By composition with the direct image of cycles by the modification $\tilde{\sigma} : \tilde{M} \to M$ we obtain, thanks to the direct image theorem for cycles by a proper map (see [B.M] ch.IV), a holomorphic map $\varphi := \tilde{\sigma}_{*}\circ \tilde{\varphi} : \tilde{N} \to \mathcal{C}_{n}^{f}(M)$.\\
 Conversely assume that $\pi : M \to N$ is a quasi-proper surjective holomorphic map  between irreducible complex spaces such that the fiber map (defined on a dense Zariski open set) is meromorphic;  let $\tau : N_{\pi} \to N$ be the modification given by the graph $N_{\pi}$ of the meromorphic fiber map of $\pi$ (see the paragraph 1.4). Then we have a holomorphic map $\tilde{\varphi} : N_{\pi} \to \mathcal{C}_{n}^{f}(M)$ which extends the fiber map. Then to prove that the strict transform $\tilde{\pi} : \tilde{M} \to N_{\pi}$ of $\pi$ by $\tau$ is  a GF map, consider the holomorphic map $\psi : N_{\pi} \to \mathcal{C}_{n}^{f}(M\times_{N}N_{\pi})$ given by $\tilde{y} \mapsto \tilde{\varphi}(\tilde{y}) \times \{\tilde{y}\}$. It is holomorphic by the product theorem for analytic families of cycles (see [B.75] or [B-M] ch. IV). As for any  $\tilde{y}$ in $N_{\pi}$ the cycle $\psi(\tilde{y})$ is in $\tilde{M} := M \times_{N, str}N_{\pi}$ this map is holomorphic with values in $\mathcal{C}_{n}^{f}(\tilde{M})$ and it is then easy to see that it is the (holomorphic) fiber map for $\tilde{\pi}$. So $\tilde{\pi}$ is a GF map.$\hfill \blacksquare$\\

\begin{defn}\label{appli. fibre}
Let $\pi : M \to N$ be a SQP map. With the notations introduced above, the map $\tilde{\varphi}  :=  \tau_{*}\circ\psi  : N_{\pi }\to \mathcal{C}_{n}^{f}(M)$ will  be called a {\bf classifying map for the family of fibers of $\pi$}.
\end{defn}

The following proposition summarizes some basic properties of a SQP map.

\begin{prop}\label{SQP basic}
Let $\pi : M \to N$ be a strongly quasi-proper holomorphic map between irreducible complex spaces. Then we have
\begin{enumerate}[i)]
\item the map $\pi$ is quasi-proper.
\item  the locus  $\Sigma_{0}$  of ``big fibers''  of $\pi$,  $\Sigma_{0} := \{y \in N \ / \  \dim \pi^{-1}(y) > n \}$, 
 where $n := \dim M - \dim N$, is a closed analytic subset in $N$. 
 \item Let $\Sigma \subset N$ be the closed analytic nowhere dense subset which is  the union of the ``big fibers'' subset and the subset of non normal points in $N$. There exists a holomorphic map $\varphi : N \setminus \Sigma \to \mathcal{C}_{n}^{f}(M)$ classifying the generic fibers of $\pi$ and  the map $\varphi$ is meromorphic along $\Sigma$.
 \item Let $\tau : N_{\pi} \to N$ the projection on $N$ of the graph of the meromorphic map $\varphi$ and $\tilde{\varphi} : N_{\pi} \to \mathcal{C}_{n}^{f}(M)$ its second projection. For each $y \in N$ we have
  $$\pi^{-1}(y) = \cup_{z \in \tau^{-1}(y)} \vert \tilde{\varphi}(z)\vert,$$
   that is to say that any ``big fiber''  of $\pi$ is filled up by limits of generic fibers.
 \end{enumerate}
 \end{prop}
 
 \parag{Proof} The properties i) and ii) has been proved in the proof of the proposition \ref{SPQ dense}. The point iii) is proved in the proposition \ref{mero et SQP}. To prove iv) consider a point $x \in M$ and consider a sequence $(x_{\nu})$ in $\pi^{-1}(N \setminus \Sigma)$ converging to $x$, where $N \setminus \Sigma$ is a  dense open set where $\varphi$ is holomorphic. This exists by irreducibility of $M$. Then $ y_{\nu} := \pi(x_{\nu}) $ converges to $y := \pi(x)$ and $\varphi(y_{\nu})$, up to pass to a sub-sequence ($\tau$ is proper), converges to a point $\tilde{y} \in \tau^{-1}(x)$. We have $x_{\nu} \in \vert \varphi(y_{\nu})\vert$ for each $\nu$ so that $x $ is in $\vert \tilde{\varphi}(\tilde{y})\vert$. $\hfill \blacksquare$

\parag{Remark} For an equidimensional  surjective holomorphic map between irreducible complex spaces to be strongly quasi-proper is equivalent to be a quasi-proper map. But if a surjective holomorphic map is not equidimensional, the condition to be quasi-proper is not strong enough; in particular it is not stable by taking the strict transform of $\pi : M \to N$ by a modification on $N$. Of course, our definition of a SQP map will be stable by such an operation. But we shall see that this notion enjoys many other stability properties. \\
The reader can also  see interesting applications of this notion in [B.13] in the construction of meromorphic quotients for a large class of non proper analytic equivalence relations.\\

\parag{Example} Let $Y := \{((a,b),(x,y)) \in \C^{2}\times  \C^{2} \ / \  a.x^{2} + b.x - a^{2}.y^{2}  = 0 \}$ and let $\pi : Y \to \C^{2}$ the map induced by the first projection. Then we have the following properties
\begin{enumerate}[I)]
\item The (algebraic) hypersurface $Y$ of $\C^{4}$ is irreducible (in fact normal and connected).
\item The map $\pi : Y \to \C^{2}$ is quasi-proper.
\item The map $\pi$ is not strongly quasi-proper. 
\item More precisely, after blowing-up the origin in $\C^{2}$ the strict transform of $\pi$ is no longer quasi-proper.
\end{enumerate}
\parag{Proof of i)} The critical set of the polynomial $P(a,b,x,y)$ if given by the following equations
\begin{equation*}
2a.x + b = 0, \quad 2a^{2}.y  = 0, \quad x^{2} - 2a.y^{2} = 0, \quad x  = 0 \tag{1}
\end{equation*}
So the subset $S := \{ a = b = x = 0\} \cup \{ x = y = b = 0\} $ which is one dimensional and is contained in $Y$ is  the singular subset of $Y$. So the singular set of $Y$ is exactly $S$. As it has codimension $2$ in $Y$, the hypersurface $Y$ is normal. We shall see that each fiber of $\pi$ is connected and then the fact that $\C^{2}\times \{0\}$ is a section of $\pi$ implies that $Y$ is connected. So $Y$ is irreducible.
\parag{Proof of ii)} First we shall describe the fibers of $\pi$. For $a.b \not= 0$   the fiber $\pi^{-1}(a,b)$ is a smooth conic containing the origin in $\C^{2}$. For  $a \not= 0$ and $b = 0$ the fiber $\pi^{-1}(a,b)$ is the union of two distinct lines through the origin. For $a = 0$ and $b \not= 0$ the fiber $\pi^{-1}(0,b)$ is the line $x  = 0$ which also contains the origin. Finally the fiber $\pi^{-1}(0,0)$ is $\C^{2}$. So each fiber is connected and contains the origin. Then the $\pi-$proper set $\C^{2}\times \{0\}$  meets every irreducible component of any fiber of $\pi$, so this map is quasi-proper.
\parag{Proof of iii)} Consider now the map $f_{s} : \C \to \C^{2}$ given by $f_{s}(b) := (b.s, b)$ where $s$ is a non zero complex number. Then for $b \not= 0$ the fiber of $\pi$ at the point $f_{s}(b)$ is given by the equation $ b.s.x^{2} + b.x - b^{2}.s^{2}.y^{2} = 0$ and assuming that $b \not= 0$ this is equivalent to 
$$ x.(s.x + 1) - b.s^{2}.y^{2} = 0 .$$
Now, for $s \not= 0$ fixed, the limit of this non degenerate conic when $b$ goes to $0$ is the union of the lines $\{x = 0\}$ and $\{x = - 1/s\}$. And it is clear that when $\vert s\vert$ goes to $0$, the component $\{x = - 1/s\}$ gets out of any compact set in $\C^{2}$. This phenomenon of escape at infinity near the point $(0,0)$ in the target of $\pi$ implies that the map $\pi$ cannot be strongly quasi-proper.
\parag{Proof of iv)} Consider now the blow-up $\tau :  X \to \C^{2}$  of the (reduced) origin in $\C^{2}$. The complex manifold $X$ is the sub-manifold
 $$X := \{((a,b),(\alpha,\beta)) \in \C^{2}\times \mathbb{P}_{1}\ / \  a.\beta = b.\alpha\}.$$
  It will be enough to show that the strict transform of $\pi$ over the chart $\{ \beta \not = 0\}$ of $X$ is not quasi-proper to achieve our goal. So let $s := \alpha/\beta$. Then we have coordinates $(s,b) \in \C^{2}$ for this chart on $X$. The total transform of $Y$ is given by the equation
$$ s.b.x^{2} + b.x - s^{2}.b^{2}.y^{2}  = 0$$
and, as the function $b$ is not generically zero on the strict transform $\tilde{Y}$ of $Y$ by $\tau$, we have
$$\tilde{Y}_{\beta \not= 0} = \{\big((s,b),(x,y)\big) \in \C^{2}\times \C^{2} \ / \  x.(s.x + 1) - b.s^{2}.y^{2} = 0 \}.$$
So the fiber of the strict transform $\tilde{\pi}$  at the point $(s,0)$ is the union of the two lines $ \{x = 0\}$ and $\{x = - 1/s\}$ for $s \not= 0$. Then it is clear that this map is not quasi-proper as an irreducible component of the fiber at $(0,s), s\not= 0$ gets out of any compact set in $\C^{2}$ when $s \not= 0$ goes to $0$.$\hfill \square$\\

\parag{Important remark} The previous example is {\bf algebraic}. And in opposition of the context where the notion of a strongly quasi-proper map has been introduced (in complex analytic geometry), this notion makes sense in algebraic geometry and is, of course, related to the behaviour of the limits of generic fibers of a map in some compactification of the map. But it is independent of the  chosen compactification and so it has to be considered also in algebraic geometry.\\
Note that the phenomenon of ``escape at infinity'' happens in the algebraic setting, but not the phenomenon of  ``infinite breaking''\footnote{We say that we have an infinite breaking when the limit in the sense of $\mathcal{C}^{loc}_{n}(M)$ of a sequence $(C_{\nu})$ of irreducible cycles has infinitely many irreducible components.}, which is purely transcendental. $\hfill \square$\\

In order to prove some more stability results for SQP maps, we shall use the following two lemma.\\

\begin{lemma}\label{exo 1}
Let $\alpha : M \to N$ and $\beta : N \to P$  be  holomorphic surjective maps between irreducible complex spaces and let  $\sigma : \tilde{P} \to P$ be a holomorphic dominant map. Let $\tilde{\beta} : \tilde{N} \to N$ the strict transform of $\beta$ by $\sigma$ and let $\tilde{\alpha} : \tilde{M} \to \tilde{N}$ the strict transform of $\alpha$ by $\tilde{\sigma}$ the projection of $\tilde{N}$ to $N$. Then  the strict transform of \  $\beta\circ\alpha$ \  by $\sigma$ is the composition of the strict transform of $\beta$ by $\sigma$ and the strict transform of $\alpha$ by $\tilde{\sigma}$.
  $$\xymatrix{\tilde{M}\ar[d]_{\tau} \ar[r]^{\tilde{\alpha}}  & \tilde{N}\ar[r]^{\tilde{\beta}} \ar[d]_{\tilde{\sigma}}& \tilde{P} \ar[d]_{\sigma} \\
M \ar[r]^{\alpha} & N \ar[r]^{\beta}& P} $$
\end{lemma}

\parag{Proof} For the ``usual'' fiber product we have a canonical isomorphism
  $$ M \times_{N}(N \times_{P}\tilde{P}) \simeq M\times_{P}\tilde{P}$$
   given by the obvious projection : the inverse is given by $(x,\tilde{z}) \mapsto (x, (\alpha(x), \sigma(\tilde{z}))$. Then our hypothesis allows to see that the strict transforms correspond via this isomorphism.$\hfill \blacksquare$\\

\begin{lemma}\label{exo 2}
Let  $M, N, \tilde{N}, P$ be irreducible complex spaces, let $\alpha : M \to N$ and  $\tau : \tilde{N} \to N$  be  holomorphic surjective maps  and $\sigma : P \to N$ a dominant holomorphic map. Note $\alpha_{1}$ and $\alpha_{2}$ the strict transforms of $\alpha$ by $\tau$ and $\sigma$ respectively. \\
 Let $P_{1} := \tilde{N}\times_{N, str}P$, and note  $\tau' : P_{1} \to \tilde{N}$ and $\sigma' : P_{1} \to P$ the projections. Then note $\alpha'_{1}$ and $\alpha'_{2}$ the strict transforms of $\alpha_{1}$ and $\alpha_{2}$ by  $\tau'$ and $\sigma'$ respectively.
$$\xymatrix{ \quad &M_{1}\ar[d]_{\alpha_{1}} \ar[r]^{\tilde{\tau}} & M \ar[d]^{\alpha} & M_{2} \ar[d]^{\alpha_{2}} \ar[l]_{\tilde{\sigma}} & \quad \\
M_{1}\times_{\tilde{N}} P_{1} \ar[ru] \ar[rrd]_{\alpha'_{1}} & \tilde{N} \ar[r]^{\tau} & N & P \ar[l]_{\sigma} & M_{2}\times_{P}P_{1} \ar[lu] \ar[lld]^{\alpha'_{2}}\\
\quad &\quad & P_{1} \ar[lu]_{\tau'} \ar[ru]^{\sigma'} & \quad & \quad } $$
Then $M_{1}\times_{\tilde{N}, str}P_{1}$ is canonically isomorphic to $M_{2}\times_{P, str}P_{1}$.
\end{lemma}

\parag{Proof} Again we have a canonical isomorphism
 $$(M \times_{N}\tilde{N})\times_{\tilde{N}} (\tilde{N}\times_{N} P_{1}) \simeq (M \times_{N}P)\times_{P}(P \times_{N} \tilde{N})$$
  given by the map $((x, \tilde{y}), (\tilde{y}, (\tilde{y}, z))) \mapsto ((x, z), (z, \tilde{y}))$;  the inverse map is given by
   $$((x,z),(z,\tilde{y})) \mapsto ((x, \tilde{y}), (\tilde{y}, (\tilde{y}, z))).$$
    Again  our hypothesis allows to see that the strict transforms correspond via this isomorphism.$\hfill \blacksquare$\\

\begin{prop}\label{strict trsf.}
 Let  $M, N, P$ be irreducible complex spaces, $\alpha : M \to N$   be a holomorphic SQP map  and $\sigma : P \to N$ a holomorphic dominant map. Then the strict transform $\alpha_{2} : M_{2} \to P$ of $\alpha$ by $\sigma$ is SQP map.
\end{prop}

\parag{Proof} By assumption, there exists a modification $\tau : \tilde{N} \to N$ such that the strict transform $\alpha_{1} : M_{1} \to \tilde{N}$ is a GF map. If $\tau' : P_{1} \to \tilde{N}$ is the strict transform of $\sigma$, it is again dominant and we can apply the lemma \ref{exo 2} to obtain that the strict transform of $\alpha_{1}$ by $\tau'$ is a GF map. Now, as the  lemma \ref{exo 2} gives that $\alpha'_{2} \simeq  \alpha'_{1}$ the conclusion follows thanks to the lemma \ref{modif. and GF}.

\begin{prop}\label{compose sqp}
 Let  $M, N, P$ be irreducible complex spaces, $\alpha : M \to N$   be a holomorphic wGF map and $\beta : N \to P$ a holomorphic SQP map. Then the composed map $\beta\circ\alpha$ is a SQP map.
 \end{prop}
 
 \parag{Proof} Let $\sigma : \tilde{P} \to P$ a modification such that the strict transform $\tilde{\beta} : \tilde{N} \to \tilde{P}$ is a GF map. Note $\tilde{\sigma} : \tilde{N} \to N$ the modification induced by the projection. Then the strict transform $\tilde{\alpha}$ of $\alpha$ by $\sigma$ is a wGF map, thanks to lemma \ref{fiber product}, and the composed map $\tilde{\beta}\circ\tilde{\alpha}$ is then a wGF map using the lemma \ref{compose} ;  then $\beta\circ \alpha$ si a SQP map by the remark following the proof of the proposition \ref{SQP et GF} and we conclude using the lemma \ref{exo 1}.$\hfill \blacksquare$

 \parag{Remarks} \begin{enumerate}
 \item Note that we dont need that the map  $\beta$ is holomorphic in the proof above :  if we compose a wGF map with a meromorphic strongly quasi-proper map\footnote{We say that a meromorphic map $ h : N --\to P$ is strongly quasi-proper if the holomorphic map $\tilde{h} : \tilde{N} \to P$ which is the projection on $P$ of the graph $\tilde{N} \subset N\times P$ of $h$ is strongly quasi-proper.} we obtain a meromorphic strongly quasi-proper map.
 \item It is not true that the composition of a modification with a SQP map is again a SQP map. The example below shows that the composition of a modification with a GF map is not quasi-proper in general.
 \end{enumerate}
 
 \parag{Example} Let $\alpha : \C^{2} \to \C$ the projection on the first coordinate. Then $\alpha$ is clearly a GF map. Consider now $\tau : X \to \C^{2}$ the blow-up in $\C^{2}$ of the reduced ideal defining the closed analytic subset $\mathbb{Z}\times \{0\}$. Then the map $\alpha\circ\tau : X \to \C$ has a fiber at $0$ which is not a finite type cycle. So this map is not quasi-proper but nevertheless  equidimensional.$\hfill \square$\\
 
We conclude this section by another important stability result for strongly quasi-proper holomorphic maps.

\begin{thm}\label{stab. tot}
Let $\alpha : M \to N$ be a SQP holomorphic map. Let $Z \subset N$ be a closed analytic irreducible subset in $N$ and $X$ an irreducible component of $\alpha^{-1}(Z)$ which is dominant on $Z$ for the map $\alpha_{\vert X}$ induced by $\alpha$. Then the map $\alpha_{\vert X} : X \to Z$ is SQP.
\end{thm}

The proof of this result will use the following proposition.

\begin{prop}\label{cycle-cycle}
Let $\alpha : M \to N$ be a holomorphic  GF map where $M$ is a pure dimensional complex space and $N$ a smooth connected manifold and define  $n := \dim M - \dim N$. Then for each integer $q \geq 0$ there is a holomorphic map $ \Phi_{q} : \mathcal{C}_{q}(N) \to \mathcal{C}_{n+q}^{f}(M)$ which is given by $\Phi_{q}(C) := (p_{2})_{*}(p_{1}^{*}(C))$, where $p_{1} : G \to N$ and $p_{2}: G \to M$ are the projections of the graph $G \subset N\times M$ of the f-analytic family of fibers of $\alpha$.
\end{prop}

\parag{Remarks}\begin{enumerate}
\item  As we know that $p_{2}$ is an isomorphism, the map $(p_{2})_{*} : \mathcal{C}_{n+q}^{f}(G) \to \mathcal{C}_{n+q}^{f}(M)$ is also  bi-holomorphic in a obvious sense (see definition \ref{hol in}). So we may identify $M$ and $G$ in the proof of the proposition.
\item As we assume that $N$ is smooth, the cycle $p_{1}^{*}(C)$ is well defined as a cycle in $G$\footnote{see [B.75] ch.VI.} as soon as the codimension of $G \cap (\vert C\vert \times M)$ in $N \times M$ is the sum of codimensions of $C$ in $N$ and of $G$ in $N\times M$, so is equal to $\dim N + \dim N - q$. This means that its dimension is equal to $\dim M - \dim N + q = n + q$. But as the fibers of $\alpha$ have pure dimension $n$ this is true for any $q-$cycle in $N$.
\end{enumerate}

\parag{Proof of the proposition \ref{cycle-cycle}} Thanks to the remark 2. above the holomorphy of the map $j\circ\Phi_{q} : \mathcal{C}_{q}(N) \to \mathcal{C}_{n+q}^{loc}(M)$, where $j : \mathcal{C}_{n+q}^{f}(M) \to \mathcal{C}_{n+q}^{loc}(M)$ is the obvious map,  is an immediate consequence of the variant for the pull-back of the intersection theorem for analytic families of cycles (see [B.75] chapter VI or [B-M 2] chapter VII). The only points to show are :
\begin{enumerate}
\item   for each compact $q-$cycle $C \in \mathcal{C}_{q}(N)$ the $(n+q)-$cycle $\Phi_{q}(C)$ is a finite type cycle;
\item the map $\Phi_{q}$ is continuous for the topology of $\mathcal{C}_{n+q}^{f}(M)$. 
\end{enumerate}
Fix $C_{0} \in \mathcal{C}_{q}(N)$ and choose a relatively compact open set $V \subset N$ such that $\vert C_{0}\vert \subset V$. As $\bar V$ is compact in $N$ there exists a relatively compact open set $W $ in $M$ such that any irreducible component of the fiber $\alpha^{-1}(y)$ for any $y \in V$ meets $W$. Take now an irreducible component $D$ of $\Phi_{q}(C')$ where $C' \in \mathcal{C}_{q}(N)$ such that  $\vert C' \vert \subset V$. We have $\vert \Phi_{q}(C')\vert = \vert p_{1}^{*}(C')\vert = G \cap (\vert C'\vert \times M) \simeq \alpha^{-1}(\vert C'\vert)$. As the restriction of $\alpha$ to $\alpha^{-1}(\vert C'\vert)$ has pure $n-$dimensional fibers over $\vert C'\vert$ its restriction to $D$ is a wGF map, thanks to the remark 2 following the definition \ref{wGF}, and $D$ is an union of irreducible components of fibers of $\alpha$ over $\vert C'\vert$. But each such irreducible component meets $W$ because $\vert C'\vert \subset V$. So $D$ meets $W$, and then $\Phi_{q}(C')$ is a finite type cycle for any such $C'$.\\
Now take a sequence $ (C_{\nu})_{\nu \geq 1}$ in $\mathcal{C}_{q}(N)$ converging to $C_{0}$. Then it is clear from the discussion above that the sequence $(\Phi_{q}(C_{\nu}))$ converges to $\Phi_{q}(C_{0})$ in $\mathcal{C}_{n+q}^{loc}(M)$ and that we have a relatively compact open set $W$ in $M$ such that any irreducible component of any $\Phi_{q}(C_{\nu})$ meets $W$. So the proposition \ref{f-conv.} implies the convergence of $(\Phi_{q}(C_{\nu}))$ to $\Phi_{q}(C_{0})$  in $ \mathcal{C}_{n+q}^{f}(M)$. This is the continuity of $\Phi_{q}$ at $C_{0}$.$\hfill \blacksquare$

\parag{Proof of the theorem \ref{stab. tot}} We begin by the description of such an $X$:\\
Let $\tau : \tilde{N} \to N$ a modification such that there exists $\varphi : \tilde{N} \to \mathcal{C}_{n}^{f}(M) $ a holomorphic map  classifying the generic fibers of  $\alpha$. Let $\tilde{Z}$ an irreducible component of $\tau^{-1}(Z)$ which is surjective onto $Z$ and let $Y := \alpha^{-1}(Z)\times_{Z,str}\tilde{Z}$. The projection $p_{2}: Y \to \tilde{Z}$ is a GF map and the projection $p_{1} : Y \to \alpha^{-1}(Z)$ is proper. Let $\tilde{X}$ be an irreducible component of $Y$ and define $X := p_{1}(\tilde{X})$. As $p_{1}$ is proper $X$ is a closed analytic irreducible subset in $\alpha^{-1}(Z)$. As $\tilde{X}$ is surjective on $\tilde{Z}$  and contains a non empty open set in $Y$, $X$ is an irreducible component of $\alpha^{-1}(Z)$ and is dominant on $Z$.\\
Conversely, if $X$ is an irreducible component of $\alpha^{-1}(Z)$ and is dominant on $Z$ then $X \times_{Z,str} \tau^{-1}(Z)$ has an irreducible component $\tilde{X}$ surjective on $X$ ($p_{1}$ and $\tau$ are proper) and $\tilde{Z} := p_{2}(\tilde{X})$ is an irreducible component of $\tau^{-1}(Z)$ which is surjective on $Z$. \\
Now fix $X, \tilde{X}$ and $\tilde{Z}$ as above. Using the geometric flattning theorem (for a proper map see [B-M] chapter IV corollary 9.3.1) and Hironaka desingularization theorem, we can find a modification $\sigma : \hat{Z} \to \tilde{Z}$ such that we have the following situation :
\begin{enumerate}[i)]
\item The proper map $\tilde{\tau} := \tau\circ \sigma : \hat{Z} \to Z$ is equidimensionnal with $\hat{Z}$ smooth an connected.
\item The strict transform $\hat{\alpha} : \hat{X} \to \hat{Z}$ of $p_{1}$ by $\sigma$ is a GF map.
\item The projection $p : \hat{X} \to X$ is proper and surjective.
\end{enumerate}
$$\xymatrix{ & & \hat{X} \ar[d] \ar[r]^{\hat{\alpha}} & \hat{Z} \ar[d]^{\sigma} &  \\  & Y\ar[d]^{p_{2}}  & \ar[l] \tilde{X} \ar[d]  \ar[r]^{p_{1}} &\tilde{Z} \ar[d] \ar[r] & \tilde{N} \ar[d]_{\tau} \\
 M & \ar[l] \alpha^{-1}(Z)& \ar[l] X \ar[r]^{\alpha_{\vert X}} & Z  \ar[r] &  N}$$
Then we have a holomorphic map $Z_{1} \to \mathcal{C}_{q}(\hat{Z})$, where $Z_{1}$ is the normalization of $Z$, classifying the generic fibers of $\tilde{\tau}$. Composing this map with the holomorphic map
 $\Psi_{q} : \mathcal{C}_{q}(\hat{Z}) \to \mathcal{C}_{n+q}^{f}(M)$ associated via the proposition \ref{cycle-cycle} to the holomorphic map $\psi : \hat{Z} \to \mathcal{C}_{n}^{f}(M)$ defined as $ \psi := \varphi_{\vert \tilde{Z}}\circ \tilde{\tau}$, gives a holomorphic map 
 $$\hat{Z} \to \mathcal{C}_{n+q}^{f}(X) \subset \mathcal{C}_{n+q}^{f}(M)$$
  which classifies the generic fibers of the restriction $ \alpha_{\vert X} : X \to Z$. This complete the proof.$\hfill \blacksquare$\\

  \section{Appendix}

  \subsection{D. Mathieu's flattening theorem.}

Our goal is to give a proof of the following result of D. Mathieu (see [M. 00]) using only the proposition \ref{SPQ dense}, its corollary \ref{invariance} and the generalization of Remmert's proper direct image theorem with values in a Banach open set (see [B-M] ch.III).\\

We shall consider a surjective holomorphic map $\pi : M \to N$ between irreducible complex spaces  which is quasi-proper (in the usual sense). So we are in the standard situation described at the begining of the paragraph 3.2. With the notations introduced there we shall assume the condition (which is our definition of a SQP map) :
\begin{itemize}
\item The projection $p : \bar \Gamma \to N$ is proper, $\hfill (@)$
\end{itemize}
where we recall that $\Gamma$ is the graph of the holomorphic map 
 $ \varphi : N \setminus \Sigma \to \mathcal{C}_{n}^{f}(M) $
classifying a f-analytic family of cycles in $M$, for $y$ generic in $N\setminus \Sigma$, and $\bar \Gamma$ its closure in $N \times  \mathcal{C}_{n}^{f}(M)$.

\begin{thm}\label{David's}
Let $\pi : M \to N$ a quasi-proper surjective holomorphic map between irreducible complex spaces and we put $n := \dim M - \dim N$. Assume that $(@)$ is satisfied. Then there exists a modification $\tau : \tilde{N} \to N$ such that the strict transform $\tilde{\pi} : \tilde{M} \to \tilde{N}$ is a geometrically f-flat map.
\end{thm}

The proof has two steps.\\ 
The first step is to prove that it is enough to prove the result locally on $\Sigma$. This will use the proper direct image theorem with values in $\mathcal{C}_{n}^{f}(M)$ which is a simple consequence (see the theorem \ref{semi-proper direct image bis}) of the generalization of Remmert's  theorem with values in a Banach space (see [B-M] ch.III). This step is given by the proposition \ref{canonical blow-up}.\\
The second step is devoted to the proof that the result is  locally true on $\Sigma$. This will be done by induction on the dimension of the fibers of $\pi$. The induction step is given by the proposition \ref{induction step}.

\begin{prop}\label{canonical blow-up}
Assume that we have a modification $\tau : \tilde{N} \to N$ such that the conclusion of the theorem is valid. Then $\bar \Gamma \subset N \times \mathcal{C}_{n}^{f}(M)$ is a closed analytic subset and the sheaf of holomorphic functions on $N \times \mathcal{C}_{n}^{f}(M)$ induces on $\bar \Gamma$  a structure of reduced complex space.
\end{prop}

\parag{Proof} Let $\tilde{\varphi} : \tilde{N} \to \mathcal{C}_{n}^{f}(\tilde{M})$ be the holomorphic map classifying the fibers of the strict transform $\tilde{\pi} : \tilde{M} \to \tilde{N}$ of $\pi$ by the modification $\tau : \tilde{N} \to N$. The composition with the direct image $\tau_{*} : \mathcal{C}_{n}^{f}(\tilde{M} \to \mathcal{C}_{n}^{f}(M)$ of cycles by the proper map $\tau$ gives a holomorphic map $\psi := \tau_{*} \circ \tilde{\varphi}$. Then define $\chi := (\tau, \psi) : \tilde{N} \to N \times \mathcal{C}_{n}^{f}(M)$. On a dense open set in $\tilde{N}$ the holomorphic map $\chi$ takes its values in $\Gamma$. So, if we prove that $\chi$ is a proper map we can conclude, first that its image is $\bar \Gamma$ and then using the proper direct image in this context we shall conclude that $\bar \Gamma$ is a (locally finite dimensional) reduced complex space with the sheaf of holomorphic functions induced from $N \times \mathcal{C}_{n}^{f}(M)$.\\
As each fiber of $\chi$ is a closed subset in a fiber of $\tau$, the fibers are compact. We have now to prove that the map $\chi$ is closed. Let $F$ be a closed set in $\tilde{N}$ and choose a sequence $(\tilde{y}_{\nu})$ be a sequence in $F$ such that the sequence $(\chi(\tilde{y}_{\nu}))$ converges to  $(y, C) \in N \times \mathcal{C}_{n}^{f}(M)$. Then, up to pass to a sub-sequence, we can assume that  the sequence $(\tilde{y}_{\nu})$ converges to some $\tilde{y} \in \tau^{-1}(y) \cap F$ by the properness of $\tau$. Then by continuity of $\chi$ we have $\chi(\tilde{y}) = (y, C)$ and so  $\chi(\tilde{y})$ is in $\chi(F)$.$\hfill \blacksquare$\\

\begin{cor}\label{univ. modif.}
Let $\pi : M \to N$ be a   holomorphic map between irreducible complex spaces, and assume that we are in the standard situation. Assume also  that for each point $\sigma$  in $\Sigma$ there exists an open neighbourhood $V_{\sigma}$ of $\sigma$ in $N$ and a modification $N_{\sigma} \to V_{\sigma}$ and a holomorphic map $\varphi_{\sigma} : N_{\sigma} \to \mathcal{C}_{n}^{f}(M)$ extending $\varphi_{\vert V_{\sigma}\setminus \Sigma}$. Then there exists a unique modification $\tau : \tilde{N} \to N$ with a holomorphic map $\tilde{\varphi} : \tilde{N} \to \mathcal{C}_{n}^{f}(M)$ extending $\varphi$ and with the following universal property :
\begin{itemize}
\item For any modification $\tau_{U} : N_{U} \to U$ of an open set $U$ in $N$ on which there exists a holomorphic map $\varphi_{U}: N_{U} \to \mathcal{C}_{n}^{f}(M)$ extending $\varphi$ on $U \setminus  \Sigma$, there exists an unique holomorphic map $\theta_{U} : N_{U} \to \tau^{-1}(U) \subset \tilde{N}$ compatible with the projections on $U$, such that $\varphi_{U} = \tilde{\varphi}\circ \theta_{U}$.
\end{itemize}
\end{cor}

\parag{Proof} This corollary is an obvious consequence of the definition of a meromorphic map of $N$  with values in $\mathcal{C}_{n}^{f}(M)$, the definition of its graph, and the universal property of the projection of this  graph on $N$ because the content of the previous proposition is precisely the meromorphy of $\varphi$ along $\Sigma$. Then the universal property gives the patching for free.$\hfill \blacksquare$ \\

Note that we only use the direct image theorem  (with values in $\mathcal{C}_{n}^{f}(M)$) only  in the proper case for the results above (and also for defining the graph of a meromorphic map with values in $\mathcal{C}_{n}^{f}(M)$).\\

\begin{prop}\label{induction step}
Assume that in the situation of the theorem \ref{David's} the map $\pi$ satisfies $(@)$. Let $y_{0}\in N$ and assume that $\dim \pi^{-1}(y_{0}) = n+k$ with $k \geq 1$, then there exists an open neighbourhood $V$ of $y_{0}$ in $N$ and a modification $\tau :  \tilde{V} \to V$ such that the strict transform $\tilde{\pi}_{V} : \tilde{M}_{V} \to  \tilde{V}$  has fibers of dimension at most $n + k - 1$, where we denote $M_{V} := \pi^{-1}(V)$ and $\pi_{V} : M_{V} \to V$ the restriction of $\pi$.
\end{prop}

First note that the following lemma has already been proved in the proof of the proposition \ref{SPQ dense}.

 \begin{lemma}\label{second stone}
 Let $\pi : M \to N$ be a holomorphic map between irreducible complex spaces. Assume that there exists a closed analytic subset $\Sigma $ with no interior point in $N$ such that the restriction of $\pi$ to $M \setminus \pi^{-1}(M) \to N\setminus \Sigma$ is quasi proper and $n-$equidimensional,  and a holomorphic map $\varphi : N \setminus \Sigma \to \mathcal{C}_{n}^{f}(\tilde{M})$ such that for $y$ generic in $N \setminus \Sigma$ we have $\varphi(y) = \vert \varphi(y)\vert = \pi^{-1}(y)$. Assume that the condition $(@)$ is satisfied. Then $\pi$ is quasi-proper.$\hfill \blacksquare$
 \end{lemma}

Remark that, thanks to the corollary \ref{invariance} we keep the hypothesis $(@)$ for $\tilde{\pi}_{V} $, and then, thanks to the lemma \ref{second stone}, the map $\tilde{\pi}$ is again quasi-proper. Remark that it is not true in general that the strict transform by a modification of a quasi-proper map is again quasi-proper (see the example following the proposition \ref{SQP basic}).

\parag{Proof of the proposition  \ref{induction step}} As $\pi$ is quasi-proper, there exists only finitely many irreducible components of $\pi^{-1}(y_{0})$. Let $\Gamma_{1}, \dots, \Gamma_{N}$ be the irreducible components of $\pi^{-1}(y_{0})$ which have dimension $n + k$. Choose on each of them a generic point $x_{i}, i \in [1,N]$, and adapted $(n+k)-$scales $E_{1}, \dots, E_{N}$, $E_{i} = (U_{i},B_{i},j_{i})$ with the following properties:
\begin{enumerate}[i)]
\item The polydisc $U_{i}$ and $B_{i}$ contain the origin ; let $pr_{i} : U_{i} \times B_{i} \to U_{i}$ be the projection.
\item  The point $x_{i}$ is in the center of $E_{I}$ and $pr_{i}(j_{i}(x_{i}) = (0,0)$.
\item $ j_{i} (\Gamma_{i}\cap j_{i}^{-1}(U_{i}\times B_{i})) = U_{i}\times \{0\}$ and $\deg_{E_{i}}(\Gamma_{i}) = 1$.
\end{enumerate}
Let $V$ be an open neighbourhood of $y_{0}$ in $N$ such that $\pi^{-1}(V) \cap j_{i}^{-1}(\bar U_{i}\times \partial B_{i}) = \emptyset$ for each $i \in [1,N]$. This exists because the compact set $\cup_{i=1}^{N} \  j_{i}^{-1}(\bar U_{i}\times \partial B_{i})$ does not meet $\pi^{-1}(y_{0})$. Let $W_{i} := \pi^{-1}(V) \cap j_{i}^{-1}(U_{i}\times B_{i})$ and consider the holomorphic map $\theta_{i}: W_{i} \to V \times U_{i}$ given by $(\pi, pr_{i}\circ j_{i})$. As $\theta_{i}^{-1}(y_{0}, 0) = \{x_{i}\}$, up to shrink $V$, we can assume that each map $\theta_{i}$ is proper and finite. We know that for generic $y$ in $V$ the analytic set $\theta_{i}(\pi^{-1}(y) \cap W_{i})$ has dimension \ $n$ \ so cannot contains $\{y\}\times U_{i}$; then the closed analytic set
$$ Z_{i} := \{ y \in V \ / \   \{y\}\times U_{i} \subset \theta_{i}(\pi^{-1}(y) \cap W_{i}) \}$$
has empty interior in $V$. Let $\mathcal{I}_{i}$ be a coherent ideal in $\mathcal{O}_{V}$ such the strict transform of $ \pi : W_{i} \to V$ by the blowing-up of $\mathcal{I}_{i}$  has no longer a fiber of dimension $\geq n + k$. See [B-M] ch.III proposition 6.1.5  for the existence of such a coherent ideal with $Supp(\mathcal{O}_{V}\big/\mathcal{I}_{i}) \subset Z_{i}$. So consider  the blow-up $\tau : \tilde{V} \to V$  in $V$ of the ideal $\mathcal{I}$ which is the product of the $\mathcal{I}_{i}$ for $i \in [1,N]$. It will give a strict transform $\tilde{\pi}_{V} : \tilde{M}_{V} \to \tilde{V}$ of $\pi_{\vert V}$  by $\tau$ whose fibers has dimension at most $n + k - 1$, concluding the proof.$\hfill \blacksquare$\\ 

\parag{Proof of the theorem \ref{David's}} First we shall consider the case where all fibers of $\pi$ have dimension at most $n$.  The lemma \ref{second stone} gives that $\pi$ is quasi-proper and so its image is closed. So the map $\pi$ is surjective and $n-$equidimensional. If $\nu : \tilde{N} \to N$ is the normalization map, then the strict transform $\tilde{\pi} : \tilde{M} \to \tilde{N}$ of $\pi$  is again $n-$equidimensional and surjective and so we have a holomorphic map $\tilde{\varphi} : \tilde{N} \to \mathcal{C}_{n}^{f}(\tilde{M})$ classifying the fibers of $\tilde{\pi}$, meaning that $\tilde{\pi}$ is a GF map.\\

Assume now that the local version of the  theorem is proved when the maximal dimension of a fiber is $n + k - 1$, with $k \geq 1$. Then the proposition \ref{induction step} shows that this local version is also true when the maximal dimension of a fiber is $n + k$. This proves the local version of the theorem. But then the corollary \ref{univ. modif.} shows that the local version of the theorem implies the theorem itself, concluding the proof.$\hfill \blacksquare$

\parag{Remark} The corollary gives a more precise result than what we state in the  theorem \ref{David's} : there exists a natural modification of $N$ which factorizes any modification of an open set $U$  in $N$ such that the map $\varphi_{U\setminus \Sigma}$ extends holomorphically on it. So the map $\varphi$ is meromorphic along $\Sigma$ and the natural modification of $N$ is simply the graph of this meromorphic map (see \ref{graph mero 2}).

\newpage

\section{Bibliography.}

 \begin{itemize}
 
 \item{[B.75]} Barlet, D. {\it Espace analytique r\'{e}duit des cycles analytiques complexes compacts d'un espace analytique complexe de dimension finie}, Fonctions de plusieurs variables complexes, II (S\'{e}m. Franois Norguet, 1974-1975), pp.1-158. Lecture Notes in Math., Vol. 482, Springer, Berlin, 1975

 \item{[B.78]} Barlet, D. \textit{Majoration du volume des fibres g\'en\'eriques et forme g\'eom\'trique du th\'eor\`{e}me d'aplatissement}, Sem. Lelong-Sloda 78-79, Lect. Notes in Math. 822, p.1-17, Springer-Verlag, Berlin.
 
  \item{[B.08]} Barlet, D. \textit{Reparam\'etrisation universelle de familles f-analytiques de cycles et f-aplatissement g\'eom\'etrique} Comment. Math. Helv. 83 (2008),\\ p. 869-888.
 
 \item{[B.13]} Barlet, D. \ {\it Quasi-proper meromorphic equivalence relations}, Math. Z. (2013), vol. 273, p. 461-484.
 
 \item{[B.15]} Barlet, D. {\it Meromorphic quotients for some holomorphic $G-$actions}. Preprint in preparation.

\item{[B-M]} Barlet, D. and Magn\'usson, J.  {\it Cycles analytiques  en g\'eom\'etrie complexe\,I: Th\'eor\`{e}me de pr\'eparation des cycles}, Cours Sp\'ecialis\'es $n^{0} 22$, Soc. Math. France (2014).

\item{[B-M 2]}  Barlet, D. and Magn\'usson, J.  {\it Cycles analytiques  en g\'eom\'etrie complexe\,II} in preparation.

\item{[H.75]} Hironaka, H. \textit{Flattening theorem in complex-analytic geometry}, Am. J. Math. 97 (1975), p. 503-547.

  \item{[K.64]} Kuhlmann, N. \textit{{\"U}ber holomorphe Abbildungen komplexer R{\"a}ume} Archiv der Math. 15 (1964), p.81-90.
 
 \item{[K.66]} Kuhlmann, N. \textit{Bemerkungen {\"u}ber  holomorphe Abbildungen komplexer R{\"a}ume} Wiss. Abh. Arbeitsgemeinschaft Nordrhein-Westfalen 33, Festschr. Ged{\"a}achtnisfeier K. Weierstrass (1966), p.495-522.
 
   \item{[M.00]} Mathieu, D. \textit{Universal reparametrization of a family of cycles : a new approach to meromorphic equivalence relations}, Ann. Inst. Fourier (Grenoble) t. 50, fasc.4 (2000) p.1155-1189.

\end{itemize}

\end{document}